\documentclass[12pt,reqno]{amsart}
\usepackage[a4paper,margin=29mm]{geometry}
\usepackage[T1]{fontenc}
\usepackage{lmodern}
\usepackage{amsmath,amssymb,amsthm,mathtools,mathrsfs}
\usepackage{microtype}
\usepackage{needspace}
\usepackage[x11names]{xcolor}
\usepackage[colorlinks=true,linkcolor=SlateBlue3,citecolor=Chartreuse4,urlcolor=Purple3]{hyperref}
\hypersetup{pdftitle={Quantitative finiteness for maximal complex totally geodesic submanifolds of non-arithmetic ball quotients},
	pdfsubject={Finite-volume non-arithmetic ball quotients and linear volume bounds},
	pdfkeywords={totally geodesic, non-arithmetic lattice, Hodge theory, toroidal compactification},
	bookmarksnumbered=true}
\renewcommand{\eqref}[1]{\textcolor{Red2}{\hyperref[#1]{(\ref*{#1})}}}
\numberwithin{equation}{section}
\newtheorem{theorem}{Theorem}[section]
\newtheorem{bigthm}{Theorem}

\newtheorem{step}{Step}

\newtheorem{proposition}[theorem]{Proposition}
\newtheorem{lemma}[theorem]{Lemma}
\newtheorem{corollary}[theorem]{Corollary}
\newtheorem{example}[theorem]{Example}
\theoremstyle{definition}
\newtheorem{definition}[theorem]{Definition}

\theoremstyle{definition}
\newtheorem{remark}[theorem]{Remark}
\newcommand{\B}{\mathbb B}
\newcommand{\C}{\mathbb C}
\newcommand{\PP}{\mathbb P}
\newcommand{\OO}{\mathcal O}
\newcommand{\KX}{K_{\overline X}}
\newcommand{\Olog}{\Omega^1_{\overline X}(\log D)}
\newcommand{\Tlog}{T^1_{\overline X}(-\log D)}
\newcommand{\Gr}{\operatorname{Gr}}
\newcommand{\GL}{\operatorname{GL}}
\newcommand{\SL}{\operatorname{SL}}
\newcommand{\SU}{\operatorname{SU}}
\newcommand{\Vol}{\operatorname{Vol}}
\newcommand{\rk}{\operatorname{rank}}
\newcommand{\Ad}{\operatorname{Ad}}
\newcommand{\tr}{\operatorname{tr}}
\newcommand{\Nash}{\operatorname{Nash}}
\newcommand{\divisor}{\operatorname{div}}
\newcommand{\Tmax}{\mathscr T^{\max}}
\newcommand{\reg}{\mathrm{reg}}

\newcommand{\specialize}{\operatorname{sp}}
\newcommand{\numeq}{\equiv_{\mathrm{num}}}
\newcommand{\cn}[1]{C_{n,#1}}
\newcommand{\degree}[1]{\deg_L(#1)}
\allowdisplaybreaks[1]
\makeatletter
\renewcommand{\@@and}{\&}
\makeatother

\title[Quantitative finiteness for totally geodesic submanifolds]{Quantitative finiteness for maximal complex totally geodesic submanifolds of non-arithmetic ball quotients}
\author{Tianzhi Hu}
\address{ School of Mathematics and Statistics, Wuhan University, Luojiashan, Wuchang, Wuhan, Hubei, 430072, P.R. China}
\email{hutianzhi@whu.edu.cn}

\author{Runze Zhang}
\address{The Institute of Mathematical Sciences and Department of Mathematics, The Chinese University of
	Hong Kong, Shatin, N.T., Hong Kong}
\email{runzezhang@cuhk.edu.hk, runze.zhang@unice.fr}
\urladdr{\href{https://sites.google.com/view/runzezhang}{https://sites.google.com/view/runzezhang}}
\date{\today}

\begin{document}
	\begin{abstract}
		We prove that the number of maximal  complex totally geodesic
		submanifolds of a smooth finite-volume non-arithmetic ball quotient is
		bounded linearly by its volume, with a constant depending only on the
		complex dimension. This gives a quantitative strengthening of the
		finiteness theorems of Bader--Fisher--Miller--Stover
		(\textit{Invent. Math.}, 2023) and Baldi--Ullmo (\textit{Ann. of Math.}, 2023).
	\end{abstract}
	\maketitle
	
	\section{Introduction}\label{sec:introduction}
	
	\subsection{Background and main result}\label{introduction 1.1}
	Let 
	$
	\B^n=\left\{z\in\C^n \,\middle|\, \sum_{j=1}^n|z_j|^2<1\right\}
	$
	be the complex unit ball, equipped with its invariant K\"ahler metric
	of holomorphic sectional curvature $-1$. Consider the ball quotient 
	$$
	X:=\Gamma\backslash\B^n,
	$$
	where $\Gamma\subset\operatorname{PU}(n,1)$ is a torsion-free lattice.
	Thus $X$ is a  complex hyperbolic manifold of finite volume.
	Recall that $\Gamma$ is \textit{arithmetic} if it is
	commensurable with
	$\phi(H(\mathbb Z)\cap H(\mathbb R)^\circ)$,
	where $H$ is a connected semisimple $\mathbb Q$-group
	and
	$\phi:H(\mathbb R)^\circ\to\operatorname{PU}(n,1)$
	is a surjective Lie group homomorphism with compact
	kernel \cite[\S12.2.4]{DeligneMostow}.
	Here $H(\mathbb Z)$ is taken in a faithful rational
	representation, and two subgroups are
	\textit{commensurable} if their intersection has
	finite index in each.
	We call $X$ \textit{arithmetic} when $\Gamma$ is,
	and \textit{non-arithmetic} otherwise.
	Non-arithmetic lattices in complex dimensions two and three were
	constructed by Mostow and Deligne--Mostow
	\cite{Mostow80,DeligneMostow}; the existence problem in higher dimensions
	remains open for $n\geq4$.
	
	Following \cite{BFMS}, a \textit{(complex) totally geodesic submanifold}
	of $X$ means a connected, positive-dimensional, properly immersed
	complex totally geodesic subspace. Such a submanifold $S\subsetneq X$
	is \textit{maximal} if it is contained in no larger proper complex
	totally geodesic submanifold of $X$. We denote the set of these maximal
	submanifolds by $\Tmax(X)$.
	The following  beautiful finiteness theorem was proved by
	Bader--Fisher--Miller--Stover and Baldi--Ullmo independently.
	
	\begin{theorem}[{\cite[Corollary 1.2]{BFMS};
			\cite[Theorem 1.2.1]{BU}}]\label{thm:known-finiteness}
		Let $n\geq2$, and let $X=\Gamma\backslash\B^n$ be a smooth
		finite-volume non-arithmetic ball quotient. Then $\Tmax(X)$ is finite.
	\end{theorem}
	
	Bader--Fisher--Miller--Stover combined superrigidity theorems with results on equidistribution
	from  homogeneous dynamics. Their theorem also covers maximal  totally real
	submanifolds of real dimension at least two; see
	\cite[Theorem 1.1]{BFMS}.
	Baldi--Ullmo gave a Hodge-theoretic proof, interpreting the finiteness
	problem in terms of unlikely intersections in a period domain for polarised $\mathbb Z$-variations of Hodge structure and
	applying the strong Ax--Schanuel theorem of Bakker--Tsimerman \cite{BT19}; see
	\cite[\S\S1.3 \& 6.1]{BU}.
	
	The analogous problem for real hyperbolic manifolds was raised by
	Reid in the mid-2000s and, independently in dimension three, by
	McMullen; see \cite[\S1]{BFMSreal}, \cite[\S8.2]{McReynoldsReid},
	\cite[Question 7.6]{DHM}.
	It was resolved by Bader--Fisher--Miller--Stover
	\cite[Theorem 1.1]{BFMSreal} and, independently for closed hyperbolic
	$3$-manifolds, by Margulis--Mohammadi
	\cite[Theorem 1.1]{MargulisMohammadi}.
	
	It is  then very natural to seek effective versions of these finiteness theorems.
	Baldi--Urbanik gave an algorithm computing all maximal complex
	totally geodesic submanifolds of non-arithmetic ball quotients whose
	lattices are given as monodromy groups of explicit smooth projective
	families. This includes the non-arithmetic Deligne--Mostow ball
	quotients \cite[Corollary 1.14 \& \S10.4]{BUr}.
	
	\medskip
	Our main result is the following \textit{quantitative} strengthening of
	Theorem \ref{thm:known-finiteness}.
	
	\begin{bigthm}\label{thm:main}
		For every integer $n\geq2$, there exists a constant $C_n>0$,
		depending only on $n$, such that every smooth finite-volume
		non-arithmetic ball quotient $X=\Gamma\backslash\B^n$ satisfies
		\begin{equation}\label{eq:main-count}
			\#\Tmax(X)\leq C_n\Vol(X).
		\end{equation}
	\end{bigthm}
	
	\medskip
	\noindent\textit{Remark.}\,\, Let $\chi(X)$ denote the topological Euler characteristic of $X$.
	The Chern--Gauss--Bonnet formula
	\cite[pp. 719--720]{HP96} gives
	$	(-1)^n\chi(X)
	=\frac{(n+1)!}{(4\pi)^n}\Vol(X)>0.
	$
	Thus Theorem \ref{thm:main} equivalently gives a linear bound
	in the \textit{signed Euler characteristic}:
	\begin{equation*}
		\#\Tmax(X)\leq\widetilde C_n\,\,(-1)^n\chi(X),
		\quad
		\widetilde C_n:=\frac{(4\pi)^n}{(n+1)!}C_n.
	\end{equation*}
	\subsection{Outline of the proof}
	
	The first main ingredient of our proof is a
	Hodge-theoretic construction.
	We pass  to a suitable finite cover, still denoted
	by $X$. Non-arithmeticity provides a nonunitary Galois
	conjugate of the standard lattice representation,
	which underlies a polarizable complex variation
	of Hodge structure \cite{BFMS,BU}.
	We use this variation to construct a nonzero global
	Wronskian section $p$ on the projective compactification
	$$
	Z:=\PP\bigl(T_{\overline X}^1(-\log D)\bigr)
	$$
	of $\PP(T_X)$, where $(\overline X,D)$ is the toroidal
	compactification of $X$.
	The section $p$ vanishes on the Zariski closure
	$\overline{A_S}\subset Z$ of the geodesic tangent incidence
	$A_S\subset Z^\circ=\PP(T_X)$ of every  totally
	geodesic submanifold $S\subsetneq X$; see Proposition \ref{prop:wronskian}.
	We remark that non-arithmeticity is essential to this
	construction; see Remark \ref{rem:Arithmetic}.
	
	Using positivity properties, we construct a very ample
	line bundle $L$ on $Z$, with $\deg_L(Z)$ bounded linearly
	by $\Vol(X)$; see Lemma \ref{lem:uniform-generation}.
	After a suitable twist, we obtain a nonzero section
	$$
	p'\in H^0(Z,L^{C_n^\prime})
	$$
	with the same vanishing property, where $C_n^\prime$ depends
	only on $n$.
	Thus all the closures $\overline{A_S}$ lie in the divisor
	$\divisor(p')$, whose $L$-degree is bounded linearly
	by $\Vol(X)$; see Corollary \ref{cor:initial-divisor}.

	\medskip
	\medskip

	The second main ingredient comes from the dynamics of
	complex geodesics. Tangent lifts of complex geodesic
	disks define the geodesic spray, which extends to
	$$
	\iota:\OO_Z(-1)\longrightarrow T_Z
	$$
	by Proposition \ref{prop:logarithmic-spray}.
	Every tangent incidence $A_S$ is invariant under
	this spray.
	Conversely, Ratner's orbit-closure theorem
	\cite{Ratner}, together with some arguments
	from \cite[Lemmata 2.7 \&  8.2]{BFMS},
	identifies every proper irreducible closed algebraic
	subvariety of $Z^\circ$ invariant under the spray
	with some $A_S$;
	see Proposition \ref{prop:invariant-incidence}.
	This characterization allows us to identify maximal
	totally geodesic submanifolds through the spray map $\iota$.
	
	\medskip
	\medskip
	
	The above two steps provide a finite descent procedure
	given in Section \ref{subsec:finitedescent} and give
	a different proof of the finiteness
	Theorem \ref{thm:known-finiteness}.
	For an irreducible subvariety $W\subset Z$
	of dimension $d$ meeting $Z^\circ$,
	the Zariski closure of its Gauss graph
	over $W_{\reg}$ records the tangent
	$d$-planes at smooth points of $W$ and their limits
	in the Grassmann bundle
	$$
	\mathsf p: \operatorname{Gr}_d(T_Z)\to Z
	$$
	(see Definition \ref{def:gauss-nash}).
	The invariance of $W\cap Z^\circ$
	under the geodesic spray can be tested on this closure
	using the normal spray $\iota^\perp$, defined by
	projecting $\mathsf p^*\iota$ onto the
	universal quotient bundle.
	Starting from the irreducible components of
	$\mathrm V(p')$, we repeatedly intersect
	the closures of the Gauss graphs
	of non-invariant subvarieties with
	$\mathrm V(\iota^\perp)$ and project
	the irreducible components of these
	intersections to $Z$ by $\mathsf p$.
	This procedure strictly decreases dimension while
	retaining the tangent incidences of all maximal
	totally geodesic submanifolds, and hence terminates
	after finitely many steps.

	\medskip
	\medskip
	Finally, we use intersection theory to control degrees
	throughout the descent in Section \ref{subsec:degree-descent}.
	We obtain an effective bound for the sum of the
	degrees of the closures of all
	maximal proper invariant subvarieties;
	see Proposition \ref{prop:degree-bound}.
	This ultimately yields the quantitative bound asserted
	in the main theorem in Section \ref{sec:volume}.

	\medskip
	\medskip
	This paper is organized  as follows.
	
	\begingroup
	\setlength{\parskip}{0pt}
	\tableofcontents
	\endgroup
	
	\section{The Hodge-theoretic Wronskian}
	\label{sec:wronskian}
	Throughout the remainder of the paper, we always assume $\textrm{dim} X=n\geq2$. 
	
	Baldi--Ullmo characterized complex totally geodesic subvarieties
	of non-arithmetic ball quotients as the positive-dimensional
	special subvarieties associated with a suitable polarised
	$\mathbb Z$-VHS \cite[Theorem 1.3.2]{BU}.
	We likewise exploit the Hodge-theoretic properties of
	totally geodesic submanifolds, but our goal here is not
	to give another such characterization.
	Instead, we construct a nonzero Wronskian section
	\eqref{eq:wronskian-section} on $Z$ (see \eqref{ZZZ})
	that vanishes on the tangent incidence $A_S$ of every
	totally geodesic submanifold $S\subsetneq X$.

	Subsection \ref{subsec:geometric-hodge-background} recalls
	the geometric and Hodge-theoretic background and defines geodesic
	tangent incidences.
	In Subsection \ref{sec:galoisconjHS}, we choose the
	conjugate representation and construct the required nef
	line subbundle.
	Subsection \ref{sec:geospray} extends the geodesic
	spray across the boundary, and
	Subsection \ref{sec:Wronski} constructs
	the Wronskian section and proves its vanishing and
	nonvanishing properties.
	\subsection{Geometric and Hodge-theoretic preliminaries}
	\label{subsec:geometric-hodge-background}
	\subsubsection{Non-arithmetic ball quotients}
	
	Let $h$ be a Hermitian form of signature $(n,1)$ on $\C^{n+1}$, and set
	\[
	G:=\operatorname{PU}(h)\cong\operatorname{PU}(n,1),
	\quad
	K:=\operatorname{Stab}_G([z_0])
	=\left\{g\in G  \,\middle|\, g[z_0]=[z_0]\right\},
	\]
	for some fixed negative line $[z_0]\in\PP(\C^{n+1})$. The negative lines
	\[
	\B^n=\left\{[z]\in\PP(\C^{n+1}) \,\middle|\, h(z,z)<0\right\}
	\]
	form the projective model of the complex ball. The action of $G$ identifies $\B^n$ with the Hermitian symmetric
	space $G/K$ of noncompact type and rank one, where
	$K\cong\operatorname{U}(n)$ is a maximal compact subgroup.
	We equip $\B^n$ with its $G$-invariant K\"ahler metric normalized
	to have holomorphic sectional curvature $-1$.
	
	A \emph{totally geodesic $k$-ball} is $\B^n\cap\PP(V)$ for a
	$(k+1)$-dimensional subspace $V\subset\C^{n+1}$ on which $h$ has
	signature $(k,1)$. Its orthogonal complement $V^{\perp_h}$ is positive definite.
	The stabilizer of the subball in $\operatorname{PU}(h)$ is the
	image of
	$\operatorname{U}(V,h|_V)\times\operatorname{U}(V^{\perp_h})$
	under the natural projection
	$\operatorname{U}(h)\to\operatorname{PU}(h)$.

	A subgroup $\Gamma\subset G$ is called a \emph{lattice} if
	it is discrete and $\Gamma\backslash G$ has finite volume. The associated ball quotient
	\[
	X:=\Gamma\backslash\B^n
	\simeq\Gamma\backslash G/K
	\]
	is a complex locally symmetric orbifold,
	and is a  manifold when $\Gamma$ is torsion-free.
	
	We use the notion of  \textit{arithmetic} lattices recalled in the beginning of the Introduction.
	We call $X$ arithmetic when $\Gamma$ is,
	and non-arithmetic otherwise.

	Throughout this paper, we always assume that $\Gamma$
	is torsion-free and focus on the non-arithmetic ball quotient; see \cite{Mostow80,DeligneMostow} for several examples.

	\subsubsection{Toroidal compactifications and logarithmic tangent bundles}
	\label{subsubsec:toroidal}
	For a noncompact finite-volume ball quotient
	$X=\Gamma\backslash\B^n$, the \emph{cusps} correspond
	to $\Gamma$-orbits of parabolic fixed points
	$\eta\in\partial\B^n$. The stabilizer
	$P_\eta=\operatorname{Stab}_{\operatorname{PU}(h)}(\eta)$
	is a maximal parabolic subgroup, and
	$\Gamma_\eta=\Gamma\cap P_\eta$ is the corresponding
	cusp subgroup. After passing to a finite cover of $X$,
	we may assume that each $\Gamma_\eta$ is contained
	in the unipotent radical of $P_\eta$, which acts
	by Heisenberg translations
	\cite[\S\S1.1--1.3]{Mok}.
	
	We recall the \emph{cusp coordinates} of \cite[\S1.1]{Mok}.
	In the Siegel domain model
	\[
	\left\{(z_1,\cdots,z_{n-1},z_n)\in\C^{n-1}\times\C\, \middle| \, 
	\operatorname{Im}z_n>\sum_{i=1}^{n-1}|z_i|^2\right\},
	\]
	the center of $\Gamma_\eta$ is generated by a
	translation $z_n\mapsto z_n+\tau$, where $\tau>0$.
	Following \cite[\S1.1, equation (3)]{Mok}, put
	\begin{equation}\label{eq:cusp-coordinates}
		(w_1,\cdots,w_{n-1},w_n)
		=\left(z_1,\cdots,z_{n-1},\exp\left(\frac{2\pi i z_n}{\tau}\right)\right).
	\end{equation}
	After taking the cusp quotient, $w_n$ is a local
	normal coordinate in a suitable trivialization,
	and the toroidal boundary is given by $w_n=0$.
	
	The \emph{smooth toroidal compactification} adds disjoint
	abelian varieties $D_i$ with anti-ample normal
	bundles, so that
	\[
	X=\overline X\setminus D,\quad
	D=\coprod_iD_i;
	\]
	see \cite[\S\S1.1--1.3]{Mok} and \cite[\S2]{BT}.
	The following lemma provides a further finite
	cover with the positivity properties needed later. Its  proof will be given in Appendix \ref{appendix}.
	
	\begin{lemma}\label{lem:deep-cover}
		After a further finite normal cover, $X$ admits a smooth
		projective toroidal compactification
		$X=\overline X\setminus D$, where
		$D=\coprod_iD_i$ is a disjoint union of abelian varieties
		with anti-ample normal bundles, and
		\begin{equation}\label{eq:deep-ample}
			\KX\ \text{ and }\ \KX-D\ \text{are ample}.
		\end{equation}
		Moreover, the cover may be chosen to factor through any prescribed
		finite cover of $X$.
	\end{lemma}
	\begin{remark}\label{rem:finitecover}
		\leavevmode
		\begin{enumerate}
			\item[(a)] Passing to a finite cover does not affect the  finiteness argument. Indeed, every maximal
			totally geodesic submanifold downstairs has only finitely
			many lifts, each of which is maximal upstairs. Thus
			finiteness on the cover implies finiteness on the original
			ball quotient.
			
			\item[(b)] Passing to a finite cover may affect effective
			estimates, since the covering degree enters the volume
			bounds. We address this issue in Section \ref{sec:volume}. \textbf{Until then, to simplify notation, we continue to
				denote the chosen finite cover by $X$}.
		\end{enumerate}
	\end{remark}

	Let $E:=\Olog$ be the logarithmic cotangent bundle. In the local coordinates
	$(w_1,\ldots,w_{n-1},w_n)$ near $D$ introduced
	in \eqref{eq:cusp-coordinates}, it is generated by
	\[
	dw_1,\cdots,\ dw_{n-1},\ \frac{dw_n}{w_n},
	\]
	while its dual $E^*=\Tlog$ is generated by\footnote{In the compact case, these logarithmic
		notations are  understood with $D=\varnothing$.}
	\begin{align}\label{eq:logvectors}
		\frac{\partial}{\partial {w_1}},\cdots,\ \frac{\partial}{\partial{w_{n-1}}},\ 
		w_n\frac{\partial}{\partial{w_n}}.
	\end{align}
	
	With the convention that projective bundles
	parametrize lines, put
	\begin{equation}\label{ZZZ}
		Z:=\PP(E^*),\quad
		\pi:Z\longrightarrow\overline X,\quad
		Z^\circ:=\pi^{-1}(X)=\PP(T_X).
	\end{equation}
	Thus $Z$ is a $(2n-1)$-dimensional projective manifold.
	The tautological subline bundle is
	$\OO_Z(-1)\subset\pi^*E^*$, and its dual
	$\OO_Z(1):=\OO_Z(-1)^*$ satisfies $\pi_*\OO_Z(1)=E$.
	
	Furthermore, a vector bundle $F$ is called \textit{ample}
	(resp.\ \textit{nef}) if $\OO_{\PP(F^*)}(1)$ is ample
	(resp.\ nef) \cite[Definition 6.1.1]{LazII}.
	In particular, $E$ is nef precisely when $\OO_Z(1)$ is nef.
	
	\subsubsection{Geodesic tangent incidences of totally geodesic submanifolds}
	Let $S\subset X$ be a totally geodesic submanifold of
	dimension $1\leq k<n$. Choose a lift
	$\widetilde S\cong\B^k\subset\B^n$ and put
	$$
	\Gamma_{\widetilde S}
	:=\left\{\gamma\in\Gamma \,\middle|\, \gamma\widetilde S=\widetilde S\right\},
	\quad
	Y_S:=\Gamma_{\widetilde S}\backslash\widetilde S.
	$$
	The natural map $\nu_S:Y_S\to S\subset X$ is the \textit{normalization
		map} and has generic degree one onto $S$. The induced complex
	hyperbolic metric on $Y_S$ is complete, and $Y_S$ has finite
	volume; see, e.g., \cite[Lemma 8.2(3)]{BFMS}.
	
	\begin{definition}
		The \emph{(geodesic) tangent incidence} $A_S\subset Z^\circ$ is defined as the image of
		$$
		\PP(T_{Y_S})\longrightarrow Z^\circ=\PP(T_X),
		\quad
		(y,[v])\longmapsto\bigl(\nu_S(y),[d\nu_S(v)]\bigr).
		$$
	\end{definition}
	
	Thus $A_S$ parametrizes the tangent lines to the local branches
	of $S$. Write $\overline{A_S}$ for its Zariski closure in $Z$.

	The following lemma describes the extension of $\nu_S$ and its
	tangent incidence map across the boundary. Its proof is also provided
	in Appendix \ref{appendix}.
	\begin{lemma}\label{lem:toroidal-incidence}
		The map $\nu_S$ extends to a finite morphism
		$f:\overline{Y_S}\to\overline X$, where
		$(\overline{Y_S},D_{Y_S})$ is the smooth projective
		toroidal compactification of $Y_S$. Its logarithmic
		differential is a subbundle injection
		$$
		d_{\log}f:T^1_{\overline{Y_S}}(-\log D_{Y_S})
		\hookrightarrow f^*\Tlog.
		$$
		Moreover, 	its projectivization
		$$
		\psi_S:\PP\bigl(T^1_{\overline{Y_S}}(-\log D_{Y_S})\bigr)
		\longrightarrow Z
		$$
		is finite and birational onto $\overline{A_S}$, with
		$$
		\psi_S^*\OO_Z(1)
		\simeq
		\OO_{\PP(T^1_{\overline{Y_S}}(-\log D_{Y_S}))}(1).
		$$
		In particular, $\overline{A_S}$ is an irreducible 
		subvariety of dimension $2k-1$, and
		$\overline{A_S}\cap Z^\circ=A_S$.
	\end{lemma}

	\subsubsection{Local systems and variation of Hodge structures}
	Let $U$ be a finite-dimensional complex vector space.
	A representation
	$$
	\rho:\pi_1(X)\longrightarrow\GL(U)
	$$
	determines a complex local system $\mathbb U$ on $X$ and an associated
	flat holomorphic bundle
	$$
	\mathcal U:=\mathbb U\otimes_\C\OO_X,
	\quad
	\nabla:\mathcal U\longrightarrow
	\mathcal U\otimes\Omega_X^1.
	$$
	
	A \textit{polarizable complex variation of Hodge structure}
	(VHS) of weight $k\in\mathbb Z$ on $\mathbb U$
	is specified by the following conditions; see e.g. 
	\cite[\S\S1.1--1.2]{Deligne87}:
	\begin{itemize}
		\item
		There is a finite smooth decomposition
		$$
		\mathcal U=\bigoplus_{p+q=k}\mathcal U^{p,q}.
		$$
		
		\item
		The associated filtration, called the \emph{Hodge filtration}
		$$
		F^a\mathcal U:=\bigoplus_{p\geq a}\mathcal U^{p,q}
		$$
		is holomorphic and satisfies \textit{Griffiths transversality}:
		$$
		\nabla F^p\mathcal U
		\subset F^{p-1}\mathcal U\otimes\Omega_X^1.
		$$
		
		\item
		There exists a nondegenerate flat Hermitian form $Q$
		for which the Hodge summands $\mathcal U^{p,q}$ are
		mutually orthogonal and $(-1)^pQ$ is positive definite
		on each nonzero $\mathcal U^{p,q}$.
		Such a form $Q$ is called a \textit{polarization}.
	\end{itemize}
	
	\subsection{Galois conjugates and a nef line subbundle}\label{sec:galoisconjHS} 
	In this subsection, we use non-arithmeticity to choose a
	nonunitary Galois conjugate of the standard lattice
	representation, corresponding to an embedding that is
	nonidentity on the adjoint trace field.
	We then use its Hodge filtration to construct a nef line
	bundle on $\overline X$ for the Wronskian construction.
	
	\begin{proposition}\label{prop:conjugate-vhs}
		The cover in Lemma \ref{lem:deep-cover} can be chosen so that
		the lattice inclusion lifts to
		$$
		\rho_0:\Gamma\longrightarrow\SU(n,1)
		\subset\SL_{n+1}(\C).
		$$
		After conjugation, its matrices lie in a number field
		$K\subset\C$. Let
		$$
		k:=\mathbb Q\bigl(\tr\Ad\rho_0(\gamma):\gamma\in\Gamma\bigr)
		$$
		be its \textit{adjoint trace field}, where
		$\Ad(g)(u):=gug^{-1}$ on $\mathfrak{sl}_{n+1}(\C)$.
		There is an embedding $\sigma:K\hookrightarrow\C$ such that
		$\rho_\sigma:=\sigma\circ\rho_0$ satisfies:
		\begin{enumerate}
			\item[\textup{(i)}]
			$\sigma|_k\ne\operatorname{id}$;
			
			\item[\textup{(ii)}]
			its image is complex Zariski dense in $\SL_{n+1}$,
			and its local monodromy along $D$ is unipotent;
			
			\item[\textup{(iii)}]
			it underlies a polarizable complex VHS and is
			\textit{nonunitary}, that is, it preserves no positive
			definite Hermitian form.
		\end{enumerate}
	\end{proposition}
	\begin{proof}
		By residual finiteness \cite[Theorem VII]{Malcev},
		a finite-index subgroup of the inverse image of
		the original lattice in $\SU(n,1)$ avoids the
		nonidentity elements of the kernel of
		$\SU(n,1)\to\operatorname{PU}(n,1)$ and projects
		injectively. This gives the lift $\rho_0$.
		The number-field construction and finite-cover
		reductions in \cite[\S7]{BFMS} arrange that
		$\rho_0$ is defined over a number field $K$
		and has unipotent boundary monodromy.
		By Lemma \ref{lem:deep-cover}, the working cover
		can be chosen inside this subgroup.
		
		Let $\mathsf G$ be the adjoint $k$-group determined
		by $\Ad\rho_0(\Gamma)$, and let
		$v_0:k\hookrightarrow\mathbb R$ be the original
		embedding. Put
		$\mathsf G_v:=\mathsf G\otimes_{k,v}\mathbb R$.
		By \cite[Theorem 1.5]{BFMS}, $k$ is totally real
		and $\Gamma$ is integral.
		The Mostow--Vinberg arithmeticity criterion
		\cite[Corollary 12.2.8]{DeligneMostow} (also recalled in
		\cite[Theorem 2.3.1]{BU}) then gives
		$$
		\Gamma\text{ is arithmetic}
		\quad\Longleftrightarrow\quad
		\mathsf G_v(\mathbb R)\text{ is compact
			for every }v\ne v_0.
		$$
		Non-arithmeticity therefore supplies $v\ne v_0$
		with $\mathsf G_v(\mathbb R)$ noncompact.
		Extend $v$ to $\sigma:K\hookrightarrow\C$.
		Then $\sigma|_k\ne\operatorname{id}$, and the
		noncompact real Zariski closure of the adjoint
		image implies that $\rho_\sigma$ is nonunitary.
		This proves \textup{(i)} and the nonunitarity
		assertion in \textup{(iii)}.
		
		For the remaining assertions, we use the
		Galois-conjugation argument in
		\cite[\S4.2, proof of Theorem 4.2.2]{BU}.
		The representation $\rho_0$ is complex Zariski
		dense in $\SL_{n+1}$
		\cite[Proposition 2.6]{Shah}, and this density
		is preserved by $\sigma$.
		For each boundary monodromy matrix $T$,
		applying $\sigma$ to $(T-I)^{n+1}=0$
		preserves unipotence. This proves \textup{(ii)}.
		
		By adjoint cohomology vanishing for lattices in $\SU(n,1)$ due to Garland--Raghunathan
		\cite[Theorem~1.10]{GarlandRaghunathan}
		(see also \cite[Theorem~4.1.2]{BU}), we have
		$$
		H^1(\Gamma,\mathfrak{sl}_{n+1}(\C))=0
		$$
		for the action $\Ad\rho_0$. The same vanishing holds
		for $\Ad\rho_\sigma$ by Galois conjugation. Thus $\rho_\sigma$ is locally rigid as an $\SL_{n+1}$-representation.
		By \cite[Lemmas~10.13 and~2.22]{Mochizuki}, it deforms within
		$\SL_{n+1}$ to a representation underlying a polarizable complex VHS.
		By the rigidity and Zariski density,  this deformation conjugate
		to $\rho_\sigma$.
		This completes \textup{(iii)}.
	\end{proof}

	Let $\mathbb V_0$ and $\mathbb V$ be the local systems
	associated with $\rho_0$ and $\rho_\sigma$, respectively.
	Set
	\begin{equation}\label{VVVVVV}
		\mathcal V:=\mathbb V\otimes_\C\OO_X,
	\end{equation}
	with its flat connection $\nabla$.
	By Proposition \ref{prop:conjugate-vhs}\textup{(iii)},
	$\mathcal V$ carries a polarizable complex VHS with
	Hodge filtration $F^\bullet\mathcal V$.
	
	Let $p_{\max}$ be the largest index for which
	$F^{p_{\max}}\mathcal V\ne0$, and put
	\begin{equation}
		\label{sWr}
		s:=\rk F^{p_{\max}}\mathcal V,
		\quad
		\mathbb W:=\bigwedge^s\mathbb V,
		\quad
		r:=\binom{n+1}{s}.
	\end{equation}
	
	\begin{remark}
		\leavevmode
		\begin{enumerate}
			\item[(a)] We have $1\leq s\leq n$ due to the    nonunitarity in
			Proposition \ref{prop:conjugate-vhs}\textup{(iii)}, otherwise the VHS would have only
			one Hodge type, making its flat polarization definite.
			\item[(b)]  Since the boundary monodromy is unipotent by
			Proposition \ref{prop:conjugate-vhs}\textup{(ii)},
			Deligne's canonical extension $\overline{\mathcal V}$
			carries a logarithmic connection with nilpotent residues
			\cite[Chapter II, Proposition 5.2]{DeligneReg}.
			Writing $j:X\hookrightarrow\overline X$, the extended
			Hodge filtration
			$$
			F^p\overline{\mathcal V}
			:=j_*F^p\mathcal V\cap\overline{\mathcal V}
			$$
			consists of subbundles with locally free successive quotients
			\cite[Theorem A]{Deng}.
		\end{enumerate}	
	\end{remark}

	We then set
	\begin{equation}\label{eq:Hodgeline}
		A:=\det F^{p_{\max}}\overline{\mathcal V},
		\quad
		\mathcal W:=\bigwedge^s\overline{\mathcal V}.
	\end{equation}
	Taking the highest exterior power of
	$F^{p_{\max}}\overline{\mathcal V}
	\hookrightarrow\overline{\mathcal V}$
	gives a line subbundle $A\hookrightarrow\mathcal W$.
	Since $\overline X$ is projective, these bundles
	and their inclusion are algebraic
	\cite[Theorems 2--3]{GAGA}.
	We also write
	$$
	\nabla:\mathcal W\longrightarrow
	\mathcal W\otimes\Omega^1_{\overline X}(\log D)
	$$
	for the induced logarithmic connection.
	
	\begin{proposition}\label{prop:hodge-determinant-nef}
		The line subbundle $A\hookrightarrow\mathcal W$ satisfies
		\begin{equation}\label{eq:A-nef}
			A\text{ is nef},
			\quad
			\det\mathcal W\simeq\OO_{\overline X}.
		\end{equation}
	\end{proposition}
	
	\begin{proof}
		Shift the Hodge types $(p,q)$ of $\mathcal V$ \eqref{VVVVVV} to
		$(p+c,q-c)$, with $c\in\mathbb Z$, adjusting the sign
		of a chosen polarization.
		Its \textit{realification}, the underlying real local system
		with its induced polarized Hodge structure, has complexification
		$\mathbb V\oplus\overline{\mathbb V}$, with conjugate summands
		of types $(q-c,p+c)$.
		For $c\gg0$, its highest filtration bundle is the original
		$F^{p_{\max}}\mathcal V$, now indexed by $p_{\max}+c$.
		
		Since $\overline X$ is projective by
		Lemma \ref{lem:deep-cover} and the boundary monodromy is
		unipotent by Proposition \ref{prop:conjugate-vhs}\textup{(ii)},
		the extension $F^{p_{\max}}\overline{\mathcal V}$ is nef
		by \cite[Corollary 1.2]{FF}.
		Its determinant $A$, defined in \eqref{eq:Hodgeline},
		is therefore nef \cite[Theorem 6.2.12]{LazII}.
		
		Since $\rho_\sigma$ takes values in $\SL_{n+1}$ by
		Proposition \ref{prop:conjugate-vhs}, we have
		$\det\mathbb V\simeq\C$, and hence
		$$
		\det\mathbb W
		\simeq(\det\mathbb V)^{\binom{n}{s-1}}
		\simeq\C.
		$$
		Tensor compatibility of the canonical extension
		\cite[Chapter II, Proposition 5.2(d)]{DeligneReg}
		gives $\det\mathcal W\simeq\OO_{\overline X}$,
		compatibly with the connections.
	\end{proof}
	
	\subsection{The tautological foliation and  geodesic spray}\label{sec:geospray}
	Recall $Z^\circ=\PP(T_X)$.
	The tautological subline bundle
	$\OO_{Z^\circ}(-1)=\OO_Z(-1)|_{Z^\circ}$ has fiber
	$$
	\OO_{Z^\circ}(-1)_\zeta
	=\C v\subset T_xX,
	\quad \zeta=(x,[v])\in Z^\circ.
	$$
	
	For each $\zeta=(x,[v])\in Z^\circ$, let
	$\Delta_\zeta$ be the unique germ of a complex geodesic
	disk through $x$ tangent to $\C v$.
	Its tangent lift is
	\[
	\widehat\Delta_\zeta
	:=\left\{(y,[T_y\Delta_\zeta])\,\,\middle |\,\, y\in\Delta_\zeta\right\}
	\subset Z^\circ.
	\]
	\begin{definition}[Tautological foliation and geodesic spray]
		\label{def:tautological-foliation-spray}
		\leavevmode
		\begin{enumerate}
			\item The holomorphic line field given, for $\zeta\in Z^\circ$, by
			$$
			\mathscr F(\zeta):=T_\zeta\widehat\Delta_\zeta
			\subset T_\zeta Z^\circ
			$$
			defines a holomorphic foliation $\mathscr F$ on $Z^\circ$,
			called the \emph{tautological foliation}
			\cite[\S3.1]{KM}; see also \cite[\S1.3]{MokGeodesic}.
			
			\item The differential of the (restricted) projection $\pi|_{Z^\circ}:Z^\circ\to X$ (see \eqref{ZZZ})  induces an
			isomorphism of line bundles
			\[
			d\pi|_{\mathscr F}:
			\mathscr F\xrightarrow{\sim}\OO_{Z^\circ}(-1),
			\]
			since, for every $\zeta=(x,[v])\in Z^\circ$,
			\[
			d\pi_\zeta:\mathscr F(\zeta)
			\xrightarrow{\sim}\C v
			=\OO_{Z^\circ}(-1)_\zeta.
			\]
			The inverse of this isomorphism, followed by the inclusion
			$\mathscr F\hookrightarrow T_{Z^\circ}$, defines
			the \emph{geodesic spray}
			\[
			\iota^\circ:\OO_{Z^\circ}(-1)
			\longrightarrow T_{Z^\circ}.
			\]
		\end{enumerate}
	\end{definition}
	
	\begin{proposition}\label{prop:logarithmic-spray}
		The spray extends to a morphism
		$$
		\iota:\OO_Z(-1)\longrightarrow T_Z,
		$$
		nowhere zero on $Z^\circ$ and tangent to $\pi^{-1}D$.
		Moreover, for every local holomorphic section $v$ of
		$\OO_Z(-1)$, the operator $\nabla_{\iota(v)}$
		preserves holomorphic local sections of $\pi^*\mathcal W$.
	\end{proposition}

	\begin{proof}
		Fix $\zeta=(x,[v])\in Z^\circ$ and choose a
		nonzero representative $v\in T_xX$.
		Locally lift $x$ and $v$ to $\B^n$, and use
		the affine Siegel coordinates $(z',z_n)$ of
		\eqref{eq:cusp-coordinates}. Write
		\[
		x=(z',z_n),\quad
		v=\sum_{j=1}^{n-1}a_j\frac{\partial}{\partial{z_j}}
		+b_0\frac{\partial}{\partial{z_n}}.
		\]
		Since the complex geodesic disk $\Delta_\zeta$
		is the intersection of $\B^n$ with a projective
		line, its lift is locally parametrized by
		\[
		t\longmapsto(z'+ta,z_n+tb_0),
		\quad a=(a_1,\ldots,a_{n-1}).
		\]
		In particular, its tangent direction $[a:b_0]$
		is constant along this parametrization.
		
		Passing to the toroidal coordinates
		$(w',w_n)$ of \eqref{eq:cusp-coordinates},
		put $b:=2\pi i b_0/\tau$. Then 
		\[
		w'(t)=w'+ta,\quad
		w_n(t)=w_n e^{bt}.
		\]
		Relative to the logarithmic tangent frame \eqref{eq:logvectors}, the tangent lift
		$\widehat\Delta_\zeta$ is locally parametrized
		in $Z^\circ$ by
		\[
		t\longmapsto
		\bigl(w'+ta,w_ne^{bt},[a:b]\bigr).
		\]
		Its tangent vector at $t=0$ is therefore\footnote{Since the fiber direction $[a:b]$ remains constant
			along the lifted curve, this tangent vector has
			no vertical component in these local
			projective-bundle coordinates.}
		\begin{equation}\label{eq:cusp-spray}
			\delta=
			\sum_{j=1}^{n-1}a_j\frac{\partial}{\partial{w_j}}
			+bw_n\frac{\partial}{\partial{w_n}}.
		\end{equation}
		By construction,
		$\delta_\zeta\in
		T_\zeta\widehat\Delta_\zeta=\mathscr F(\zeta)$,
		and $d\pi_\zeta(\delta_\zeta)=v$.
		Hence $\delta_\zeta=\iota^\circ_\zeta(v)$
		by the definition of the geodesic spray; see Definition \ref{def:tautological-foliation-spray}(2).
		
		We now extend $\iota^\circ$ across the boundary using
		\eqref{eq:cusp-spray}.
		The expression is holomorphic at $w_n=0$ and linear
		in the homogeneous fiber coordinates $(a,b)$.
		After normalizing a nonzero fiber coordinate, it defines
		a local bundle morphism
		$$
		\OO_Z(-1)\longrightarrow T_Z.
		$$
		Since these local morphisms restrict to $\iota^\circ$
		on the dense interior, they agree on overlaps and
		glue to a morphism $\iota$ on $Z$.
		
		Moreover, $\delta(w_n)=bw_n$, so $\iota$ is
		tangent to $\pi^{-1}D$. Its restriction to
		$Z^\circ$ is nonzero, since
		$d\pi\circ\iota^\circ$ is the 
		inclusion.
		The morphism $\iota$ is algebraic since $Z$ is projective.
		
		Finally, in a local frame of
		$\mathcal W$, the logarithmic connection
		has the form
		\[
		\nabla=d+
		\sum_{j=1}^{n-1}\Theta_j\,dw_j
		+\Theta_n\frac{dw_n}{w_n},
		\]
		where the $\Theta_j$ are holomorphic matrix-valued functions.  By \eqref{eq:cusp-spray}, we have
		\[
		\nabla_\delta
		=\delta+\sum_{j=1}^{n-1}a_j\Theta_j
		+b\Theta_n,
		\]
		which has holomorphic coefficients across
		$w_n=0$. Thus $\nabla_\delta$ preserves
		holomorphic local sections, as do its
		successive iterates. 
	\end{proof}

	\subsection{The Wronskian section}\label{sec:Wronski}
	Recall the line subbundle $A\hookrightarrow\mathcal W$
	from \eqref{eq:Hodgeline}, with
	$r=\rk\mathcal W$ as in \eqref{sWr}.
	Put $\nu:=r(r-1)/2$.
	For local frames $a$ of $\pi^*A$ and $v$ of $\OO_Z(-1)$,
	write $\nabla_v:=\nabla_{\iota(v)}$, where $\nabla$
	is the pulled-back logarithmic connection and $\iota$
	is the spray of Proposition \ref{prop:logarithmic-spray}.
	The associated \textit{local Wronskian expression} is
	\begin{equation}\label{eq:wronskian-local}
		p(a,v):=a\wedge\nabla_va\wedge\cdots
		\wedge\nabla_v^{r-1}a.
	\end{equation}

	Choose a flat generator of $\det\mathcal W$, identifying $p(a,v)$
	with a local holomorphic function. On a leaf, choose $t$ with
	$\iota(v)=\frac{\partial}{\partial t}$ and a frame satisfying $\nabla e_i=0$ and
	$e_1\wedge\cdots\wedge e_r=1$. If $a=\sum_i f_i(t)e_i$, then
	$$
	p(a,v)=\det\bigl(f_i^{(j-1)}(t)\bigr)_{1\leq i,j\leq r},
	$$
	the classical Wronskian of $f_1,\ldots,f_r$; see e.g. 
	\cite[p. 722]{BostanDumas}.
	
	\begin{proposition}\label{prop:wronskian}
		The local expressions \eqref{eq:wronskian-local}
		define a nonzero global algebraic section,
		called the \textit{(global) Wronskian section},
		\begin{equation}\label{eq:wronskian-section}
			p\in H^0\bigl(Z,\OO_Z(\nu)\otimes\pi^*A^{-r}\bigr),
		\end{equation}
		vanishing on $\overline{A_S}$ for every totally geodesic submanifold $S\subsetneq X$.
	\end{proposition}
	
	\begin{proof}	The proof  will be divided into three steps.
		\setcounter{step}{0}
		\begin{step}{Globalness.}
		\end{step}
		Using the flat trivialization of $\det\mathcal W$ from
		the proof of Proposition \ref{prop:hodge-determinant-nef},
		the Leibniz rule gives $p(fa,hv)=f^rh^\nu p(a,v)$
		for local holomorphic units $f,h$.
		By Proposition \ref{prop:logarithmic-spray},
		these expressions are holomorphic across $\pi^{-1}D$
		and hence define a global holomorphic section
		\eqref{eq:wronskian-section}.
		Since $Z$ is projective, this section is algebraic
		by \cite[Theorem 1]{GAGA}.

		\begin{step}{Vanishing on $\overline{A_S}$.}\label{step 2222}
		\end{step}
		Let $k:=\dim S<n$ and choose a lift
		$\widetilde S=\B^n\cap\PP(P_0)$, where
		$P_0\subset\C^{n+1}$ is a $(k+1)$-dimensional
		subspace on which $h$ has signature $(k,1)$.
		The spaces $P_0$ and $P_0^{\perp_h}$ define local systems
		$\mathbb P_0,\mathbb N_0$. Writing $|_{Y_S}$ for pullback
		by $\nu_S$, we have
		$$
		\mathbb V_0|_{Y_S}=\mathbb P_0\oplus\mathbb N_0,
		\quad
		\rk\mathbb P_0=k+1,\quad \rk\mathbb N_0=n-k.
		$$
		By \cite[Proposition 2.6]{Shah}, the real Zariski closure
		of $\Gamma_{\widetilde S}$ contains the special unitary
		group of $P_0$, acting trivially on $P_0^{\perp_h}$.
		Thus $\mathbb P_0$ is irreducible and nonunitary, while
		$\mathbb N_0$ is unitary.
		
		Let $B_0$ and $B$ be the algebras of global flat endomorphisms
		of $\mathbb V_0|_{Y_S}$ and $\mathbb V|_{Y_S}$, respectively.
		Since $Y_S$ admits the compactification of
		Lemma~\ref{lem:toroidal-incidence}, the theorem of the fixed part
		\cite[Theorem~7.22]{Schmid}, in its complex form
		\cite[p.~2]{BakkerMatsushita}, gives
		$B=\bigoplus_a B^{a,-a}$.
		
		The commutant $B_0$ is defined over the number field $K$
		of Proposition \ref{prop:conjugate-vhs}.
		Since $\mathbb P_0$ is irreducible and nonunitary whereas
		$\mathbb N_0$ is unitary, the projector $\varpi$ onto
		$\mathbb P_0$ is central.
		The central idempotents form a finite algebraic set over $K$.
		Extending $\sigma$ to $\overline K$ gives a central
		idempotent $\varpi^\sigma\in B$ and a flat splitting
		$\mathbb V|_{Y_S}=\mathbb P\oplus\mathbb N$
		of the same ranks.
		
		The Hodge grading is compatible with composition and defines
		a $\C^*$-action on $B$ by algebra automorphisms,
		$b\mapsto t^ab$ for $b\in B^{a,-a}$.
		Connectedness forces this action to fix every central
		idempotent. Hence $\varpi^\sigma\in B^{0,0}$
		and preserves $F^\bullet$.
		
		For the flat bundles $\mathcal P,\mathcal N$ of
		$\mathbb P,\mathbb N$, put
		$i:=\rk(F^{p_{\max}}\mathcal V|_{Y_S}\cap\mathcal P)$.
		Since $\varpi^\sigma$ preserves the Hodge filtration,
		$$
		F^{p_{\max}}\mathcal V|_{Y_S}
		=
		\bigl(F^{p_{\max}}\mathcal V|_{Y_S}\cap\mathcal P\bigr)
		\oplus
		\bigl(F^{p_{\max}}\mathcal V|_{Y_S}\cap\mathcal N\bigr).
		$$
		Taking determinants and using \eqref{eq:Hodgeline} gives
		$$
		A|_{Y_S}\subset
		\bigwedge^i\mathcal P\otimes\bigwedge^{s-i}\mathcal N
		\subset\mathcal W|_{Y_S}.
		$$
		This summand is proper: since $1\leq s\leq n$ and
		$1\leq k<n$, the identity
		$$
		\binom{n+1}{s}
		=\sum_j\binom{k+1}{j}\binom{n-k}{s-j}
		$$
		has at least two positive terms.
		Since $S$ is totally geodesic, all derivatives in
		\eqref{eq:wronskian-local} along $A_S$ remain in this
		proper flat summand.
		Thus $p|_{A_S}=0$, and regularity on $Z$ gives
		$p|_{\overline{A_S}}=0$.
		
		\begin{step}{Nonvanishing.}
		\end{step}
		Suppose $p\equiv0$. We derive a contradiction to the
		irreducibility of $W:=\bigwedge^s\C^{n+1}$.
		Pull $A|_X\subset\mathcal W|_X$ back to $\B^n$.
		A flat trivialization gives $\Phi:\B^n\to\PP(W)$.
		For a projective line $\ell$ meeting the ball, set
		$U_\ell:=\sum_{x\in\ell\cap\B^n}\Phi(x)\subset W$,
		viewing each $\Phi(x)$ as a line in $W$.
		The Wronskian criterion \cite[Theorem 1]{BostanDumas}
		gives $b:=\max_\ell\dim_\C U_\ell<r$.
		On the domain $\mathcal L\subset Q:=\Gr(2,\C^{n+1})$
		of such lines, $\ell\mapsto U_\ell$ defines a meromorphic map
		$$
		\Psi:\mathcal L\dashrightarrow\Gr(b,W).
		$$
		Locally choose a point and a coordinate $t$ on each disk,
		varying holomorphically with $\ell$, and a lift
		$F(t,\ell)$ of $\Phi$.
		A nonzero wedge of $b$ derivatives
		$\partial_t^{j_1}F,\ldots,\partial_t^{j_b}F$ at that point
		gives Pl\"ucker coordinates for $U_\ell$.
		Wherever $\dim U_\ell=b$, these local maps agree.
		
		For a negative line $\Lambda\subset\C^{n+1}$, the planes
		containing it form $C_\Lambda\simeq\PP^{n-1}\subset\mathcal L$.
		The domain $\mathcal L$ is connected, being the image of
		$\PP(T_{\B^n})$.
		Since $Q$ is simply connected with Picard rank one,
		every nonzero effective divisor meets $C_\Lambda$.
		Applying the implication $(13)\Rightarrow(7)$ of
		\cite[Theorem 2 \& Comments 3]{KollarNeighborhoods}
		to $C_\Lambda\subset\mathcal L\subset Q$, one obtains
		${\mathcal{M}er}(\mathcal L)=\mathcal{R}{at}(Q)$,
		where these denote the fields of meromorphic functions
		on $\mathcal L$ and rational functions on $Q$, respectively.
		Applying this to the Pl\"ucker coordinate ratios extends
		$\Psi$ to a rational map on $Q$.
		
		\textit{\underline{Claim.}}
		\begin{equation}\label{eq:product-monodromy}
			\overline{(\rho_0,\rho_\sigma)(\Gamma)}^{\mathrm{Zar}}
			=\SL_{n+1}\times\SL_{n+1}.
		\end{equation}
		
		\begin{proof}[Proof of Claim]
			Both projections are surjective.
			If the closure were proper, Goursat's lemma,
			as used in \cite[proof of Theorem 5.1.4]{BU},
			would identify the two projective representations
			up to an automorphism of $\operatorname{PSL}_{n+1}$. 
			Such automorphisms preserve adjoint traces, giving
			$$
			\sigma\bigl(\tr\Ad\rho_0(\gamma)\bigr)
			=\tr\Ad\rho_0(\gamma)
			\quad(\gamma\in\Gamma).
			$$
			Since these traces generate the adjoint trace field,
			this contradicts
			Proposition \ref{prop:conjugate-vhs}\textup{(i)}.
		\end{proof}
		Since $\Phi(\gamma x)
		=(\bigwedge^s\rho_\sigma(\gamma))\Phi(x)$,
		the Zariski closure of the graph of $\Psi$ is invariant
		under $(\rho_0,\rho_\sigma)(\Gamma)$.
		Its stabilizer is Zariski closed, so
		\eqref{eq:product-monodromy} makes it invariant under
		$\SL_{n+1}\times\SL_{n+1}$.
		The second factor fixes $Q$ and therefore preserves
		the $b$-plane $\Psi(\ell)$ for general $\ell\in Q$.
		Since $0<b<r$, this contradicts the irreducibility of $W$.
		Thus $p\not\equiv0$.
	\end{proof}

	\begin{remark}\label{rem:Arithmetic}
		The non-arithmeticity assumption is essential to
		our construction. Indeed, by the arithmeticity criterion
		recalled in the proof of
		Proposition \ref{prop:conjugate-vhs}, if $\Gamma$ is
		arithmetic, then $\mathsf G_v(\mathbb R)$ is
		compact at every real place $v\neq v_0$ of
		the adjoint trace field $k$.
		Consequently, every Galois conjugate
		$\rho_\sigma$ with $\sigma|_k\neq\operatorname{id}$
		is unitary. Thus Proposition \ref{prop:conjugate-vhs} fails and the nonunitary
		Galois conjugate required in our construction
		does not exist in the arithmetic case.
	\end{remark}

	\section{Positivity and a common equation}
	\label{sec:positivity}
	Proposition \ref{prop:wronskian} gives a nonzero Wronskian section $p$
	vanishing on $\overline{A_S}$ for every proper totally geodesic submanifold $S\subset X$. In this section, 
	we will select a very ample line bundle $L$ on $Z$ \eqref{ZZZ}  and use
	the nefness of $A$ from Proposition \ref{prop:hodge-determinant-nef} to obtain a nonzero section
	$p'\in H^0(Z,L^{\cn{5}})$ with the same vanishing property. Here,
	the exponent $\cn{5}$ depends only on $n$, and the degree
	of $\divisor(p')$ is bounded linearly by $\Vol(X)$.
	The construction also gives global generation of
	$T_Z\otimes L^{\cn{3}}$, which will supply the
	tangency equations in Section \ref{sec:descent}.

	For a line bundle $L$ on a projective variety $Y$
	and a $d$-cycle $\alpha=\sum_j n_j[W_j]$, define
	$$
	\degree{\alpha}
	:=\sum_j n_j\int_{W_j}c_1(L)^d,
	$$
	where the $W_j$ are distinct irreducible
	$d$-dimensional subvarieties and $n_j\in\mathbb Z$.
	The cycle $\alpha$ is called \textit{effective}
	if all $n_j\geq0$.
	For an irreducible subvariety $W\subset Y$, write
	$\degree{W}:=\degree{[W]}$.
	
	\begin{proposition}\label{prop:toroidal-properties}
		The bundle $E=\Omega_{\overline X}^1(\log D)$ is nef. Moreover,
		\begin{equation}\label{eq:canonical-volume}
			\int_{\overline X}c_1\!\left(\OO_{\overline X}(\KX+D)\right)^n=n!\left(\frac{n+1}{4\pi}\right)^n\Vol(X).
		\end{equation}
	\end{proposition}
	
	\begin{proof}
		Nefness of $E$, equivalently of $\OO_Z(1)$,
		follows from \cite[Theorem 3]{Cadorel}.
		Let $\omega$ be the hyperbolic K\"ahler form.
		The induced metric on $K_X$ has first Chern form
		$(n+1)\omega/(4\pi)$.
		The logarithmic Chern--Weil formula
		\cite[Proposition 2.2.2]{Cadorel},
		based on \cite[Theorem 1.4]{Mumford},  then gives
		\eqref{eq:canonical-volume}.
	\end{proof}
	
	\begin{lemma}\label{lem:uniform-generation}
		There are integers $\cn{2},\cn{3}>0$, depending only on $n$, such that
		$L:=\OO_Z(1)\otimes\pi^*\KX^{\cn{2}}$ is very ample and both
		$\pi^*\KX^{\cn{2}}$ and $T_Z\otimes L^{\cn{3}}$ are globally generated.
		Moreover,
		\begin{equation}\label{eq:ambient-degree}
			\begin{aligned}
				\deg_L(Z)
				&\leq\cn{4}\int_{\overline X}c_1\!\left(\OO_{\overline X}(\KX+D)\right)^n,\\[0.5em]
				\cn{4}
				&:=\sum_{j=0}^n\binom{2n-1}{n-1+j}\cn{2}^{\,n-j}.
			\end{aligned}
		\end{equation}
	\end{lemma}
	
	\begin{proof}
		Set $$\cn{1}:=12n^n, \quad \cn{2}:=(n+2)\cn{1}\quad
		\text{and}\quad  \cn{3}:=6n.$$
		Since $\KX$ is ample by Lemma \ref{lem:deep-cover},
		the bundle $B:=\KX^{\cn{1}}$ is very ample by
		\cite[Corollary 12.10]{Demailly}.
		For a nef bundle $F$ and an ample line bundle $B'$
		on a  projective manifold $Y$, Griffiths vanishing gives
		\begin{equation}\label{eq:griffiths-vanishing}
			H^q(Y,K_Y\otimes F\otimes\det F\otimes B')=0
			\quad(q>0)
		\end{equation}
		\cite[Variant 7.3.2]{LazII}.
		We use the following regularity criterion for a vector
		bundle $G$ and a very ample line bundle $H$ on $Y$
		\cite[Theorem 1.8.5]{LazI}:
		\begin{equation}\label{eq:regularity}
			H^j(Y,G\otimes H^{-j})=0\quad(j>0)
			\quad\Longrightarrow\quad
			G\text{ is globally generated}.
		\end{equation}
		
		We now complete the proof in the following three steps.
		\setcounter{step}{0}
		\begin{step}{Very ampleness.}
		\end{step}\nopagebreak[4]
		Apply \eqref{eq:griffiths-vanishing} to $E$ with
		$B':=\KX^{t-3}\otimes\OO_{\overline X}(\KX-D)$.
		By \eqref{eq:deep-ample},
		\begin{equation}\label{eq:cotangent-vanishing}
			H^q(\overline X,E\otimes\KX^t)=0
			\quad(q>0,\ t\geq4).
		\end{equation}
		Since $(n+1-j)\cn{1}\geq4$ for $1\leq j\leq n$,
		\eqref{eq:cotangent-vanishing} and \eqref{eq:regularity}
		give that $E\otimes B^{n+1}$ is globally generated.
		Its evaluation map embeds $Z$ into
		$\overline X\times\PP^N$.
		The second projection pulls $\OO_{\PP^N}(1)$ back to
		$\OO_Z(1)\otimes\pi^*B^{n+1}$.
		Combining with the $B$-embedding of $\overline X$
		and the Segre embedding proves that
		$L=\OO_Z(1)\otimes\pi^*B^{n+2}$ is very ample.
		Moreover, $\pi^*\KX^{\cn{2}}=\pi^*B^{n+2}$
		is globally generated, since $B$ is very ample.
		
		\begin{step}{Tangent generation.}
		\end{step}\nopagebreak[4]
		By \eqref{eq:deep-ample} and
		\cite[Theorem 6.2.12]{LazII},
		$$
		E^*\otimes\KX^2
		\simeq\bigwedge^{n-1}E
		\otimes\OO_{\overline X}(\KX-D)
		\quad\text{is nef}.
		$$
		The ordinary twist $T_{\overline X}\otimes\KX^2$
		is also nef.
		For $f:C\to\overline X$ from a smooth projective curve
		with $f(C)\not\subset D$, the induced map from
		$f^*(E^*\otimes\KX^2)$ to any quotient line $Q$ of
		$f^*(T_{\overline X}\otimes\KX^2)$ has image $Q(-P)$
		for an effective divisor $P$.
		Nefness gives $\deg Q\geq\deg Q(-P)\geq0$
		\cite[Proposition 6.1.18]{LazII}.
		On $D_i$, use
		$$
		0\longrightarrow
		(\KX^2|_{D_i})^{\oplus(n-1)}
		\longrightarrow
		(T_{\overline X}\otimes\KX^2)|_{D_i}
		\longrightarrow
		\KX|_{D_i}
		\longrightarrow0,
		$$
		where $\KX|_{D_i}=\OO_{D_i}(-D_i)$ is ample.
		
		The relative Euler and tangent sequences
		\begin{equation}\label{Euler sequence}
			0 \longrightarrow\OO_Z
			\longrightarrow\pi^*E^*\otimes\OO_Z(1)
			\longrightarrow T_{Z/\overline X}
			\longrightarrow0
		\end{equation}
		and
		\begin{equation}\label{relative tangent}
			0 \longrightarrow T_{Z/\overline X}
			\longrightarrow T_Z
			\longrightarrow\pi^*T_{\overline X}
			\longrightarrow0,
		\end{equation}
		show that $F:=T_Z\otimes\pi^*\KX^2$ is nef.
		Since $\rk T_Z=2n-1$,
		$$
		K_Z\otimes F\otimes\det F
		=T_Z\otimes\pi^*\KX^{4n}.
		$$
		Since $E$ is nef by Proposition \ref{prop:toroidal-properties}
		and $\KX$ is ample by Lemma \ref{lem:deep-cover},
		the bundle
		$$
		L^t\otimes\pi^*\KX^{-4n}
		=\OO_Z(t)\otimes\pi^*\KX^{\cn{2}t-4n}
		$$
		is ample for $t>0$ and $\cn{2}t>4n$
		\cite[Proposition 6.2.11]{LazII}.
		Thus \eqref{eq:griffiths-vanishing} gives
		\begin{equation}\label{eq:tangent-vanishing}
			H^q(Z,T_Z\otimes L^t)=0
			\quad(q>0,\ t>0,\ \cn{2}t>4n).
		\end{equation}
		Taking $t=\cn{3}-j$, $1\leq j\leq2n-1$, in
		\eqref{eq:tangent-vanishing}, criterion
		\eqref{eq:regularity} gives that
		$T_Z\otimes L^{\cn{3}}$ is globally generated.
		
		\begin{step}{Degree.}
		\end{step}\nopagebreak[4]
		For the Chern roots $x_1,\ldots,x_n$ of $E$, put
		$$
		h_j(E):=
		\sum_{a_1+\cdots+a_n=j}
		x_1^{a_1}\cdots x_n^{a_n},
		\quad h_0(E):=1.
		$$
		The projective-bundle formula gives
		$\pi_*c_1(\OO_Z(1))^{n-1+j}=h_j(E)$
		\cite[\S\S3.1--3.3]{Fulton}.
		The difference $c_1(E)^j-h_j(E)$ is a nonnegative
		sum of Schur classes.
		Schur positivity \cite[Theorem I]{FL} extends to nef $E$
		by replacing $x_i$ with $x_i+\delta c_1(\KX)$
		and letting the positive rational number $\delta$
		tend to zero.
		It gives
		$$
		\begin{aligned}
			0&\leq
			\int_{\overline X}h_j(E)c_1(\KX)^{n-j}\\[0.5em]
			&\leq
			\int_{\overline X}c_1\!\left(\OO_{\overline X}(\KX+D)\right)^jc_1(\KX)^{n-j}\\[0.5em]
			&\leq\int_{\overline X}c_1\!\left(\OO_{\overline X}(\KX+D)\right)^n.
		\end{aligned}
		$$
		The last inequality uses the effectivity of $D$
		and the nefness of $\KX$ and $\KX+D$.
		Consequently,
		$$
		\begin{aligned}
			\deg_L(Z)
			&=\sum_{j=0}^n
			\binom{2n-1}{n-1+j}\cn{2}^{\,n-j}
			\int_{\overline X}h_j(E)c_1(\KX)^{n-j}\\[0.6em]
			&\leq\cn{4}\int_{\overline X}c_1\!\left(\OO_{\overline X}(\KX+D)\right)^n.
		\end{aligned}
		$$
	\end{proof}

	We now use the nefness of $A$ in \eqref{eq:A-nef}
	to turn $p$ into a section of some power of $L$.
	The bound \eqref{eq:ambient-degree} then controls
	the degree of its zero divisor.
	
	\begin{corollary}\label{cor:initial-divisor}
		Set
		\begin{equation}\label{eq:initial-exponent}
			\cn{5}:=2n+3+\max_{1\leq u\leq n}
			\binom{\binom{n+1}{u}}{2}.
		\end{equation}
		There is $$
		0\ne p'\in H^0(Z,L^{\cn{5}})
		$$
		vanishing on every
		$\overline{A_S}$ with $S$ proper. In particular,
		$$
		\begin{aligned}
			\deg_L\divisor(p')
			&=\cn{5}\deg_L(Z)\\[0em]
			&\leq
			\cn{5}\cn{4}\,n!
			\left(\frac{n+1}{4\pi}\right)^n\Vol(X).
		\end{aligned}
		$$
	\end{corollary}
	
	\begin{proof}
		Recall that $r=\rk\mathcal W=\binom{n+1}{s}$ as in
		\eqref{sWr}, and $\nu=r(r-1)/2$.
		By \eqref{eq:wronskian-section},
		$p\in H^0(Z,\OO_Z(\nu)\otimes\pi^*A^{-r})$.
		Put
		$$
		B:=\OO_Z(-\nu)\otimes\pi^*A^r\otimes L^{\cn{5}}.
		$$
		Taking determinants in the relative Euler sequence
		\eqref{Euler sequence} and the relative tangent sequence \eqref{relative tangent},
		and using $\det E\simeq\OO_{\overline X}(\KX+D)$,
		we obtain
		$$
		K_Z\simeq\OO_Z(-n)\otimes
		\pi^*\OO_{\overline X}(2\KX+D).
		$$
		For $1\leq j\leq2n-1$, we have
		$$
		\begin{aligned}
			B\otimes L^{-j}\otimes K_Z^{-1}
			={}&\OO_Z(\cn{5}-j-\nu+n)\\
			&\otimes\pi^*\!\left(\KX^{\cn{2}(\cn{5}-j)-3}
			\otimes\OO_{\overline X}(\KX-D)\otimes A^r\right).
		\end{aligned}
		$$
		By \eqref{eq:initial-exponent},
		$\cn{5}-j\geq\nu+4$.
		After tensoring the right-hand side with $L^{-1}$,
		the powers of $\OO_Z(1)$ and $\KX$ remain nonnegative.
		Thus $B\otimes L^{-j-1}\otimes K_Z^{-1}$
		is nef by Proposition \ref{prop:toroidal-properties},
		\eqref{eq:deep-ample}, and \eqref{eq:A-nef}.
		Tensoring with the very ample bundle $L$ of
		Lemma \ref{lem:uniform-generation} shows that
		$B\otimes L^{-j}\otimes K_Z^{-1}$ is ample.
		Kodaira vanishing gives
		$$
		H^j(Z,B\otimes L^{-j})=0
		\quad(1\leq j\leq2n-1).
		$$
		Criterion \eqref{eq:regularity} then shows that
		$B$ is globally generated.

		Choose $z\in Z$ with $p(z)\ne0$ and
		$s_0\in H^0(Z,B)$ with $s_0(z)\ne0$.
		Then $$p':=ps_0\in H^0(Z,L^{\cn{5}})$$ is nonzero
		and vanishes on every $\overline{A_S}$ by
		Proposition \ref{prop:wronskian}.
		Since $\OO_Z(\divisor(p'))\simeq L^{\cn{5}}$,
		we have
		$$
		\deg_L\divisor(p')=\cn{5}\deg_L(Z).
		$$
		The asserted bound follows from
		\eqref{eq:ambient-degree} and
		\eqref{eq:canonical-volume}.
	\end{proof}
	
	\textbf{We will fix such a section $p'$ for the remainder of the paper.}

	\section{Subvarieties invariant under the geodesic spray}
	\label{sec:invariant-subvarieties}
	In this section, we characterize the proper irreducible
	closed algebraic subvarieties of $Z^\circ$ invariant
	under the geodesic spray of
	Definition \ref{def:tautological-foliation-spray}.
	This characterization will be used in
	Section \ref{sec:descent} to identify the invariant
	subvarieties obtained by successive tangency conditions.
	
	\begin{definition}
		We call a reduced algebraic subvariety $I\subset Z^\circ$
		\textit{invariant under the geodesic spray} if
		$\iota_x(\OO_Z(-1)_x)\subset T_xI$
		for every $x\in I_{\reg}$,
		where $T_xI$ denotes the complex tangent space of $I$ at $x$.
		Equivalently, local regular vector fields generating
		the spray preserve the ideal sheaf of $I$ \cite[\S3.5]{Druel}.
	\end{definition}
	
	By Lemma \ref{lem:toroidal-incidence}, each $A_S$ is
	a closed irreducible algebraic subvariety of $Z^\circ$.
	It is invariant under the geodesic spray,
	since $S$ is totally geodesic.
	The following proposition gives the converse.
	
	\begin{proposition}\label{prop:invariant-incidence}
		Every proper irreducible closed algebraic subvariety
		of $Z^\circ$ invariant under the geodesic spray
		is $A_S$ for some totally geodesic submanifold
		$S\subsetneq X$.
	\end{proposition}
	
	\begin{proof}
		\setcounter{step}{0} 	The proof  will be divided into four steps.
		\begin{step}{Locally symmetric space.}
		\end{step}
		Let $I\subsetneq Z^\circ$ be irreducible and invariant.
		Fix $(o,\ell_0)\in\PP(T_{\B^n})$ and put $G:=\SU(n,1)$,
		$K:=\operatorname{Stab}_G(o)$ and
		$$
		K_0:=\operatorname{Stab}_G(o,\ell_0)
		=S\bigl(U(1)\times U(1)\times U(n-1)\bigr),
		\quad Z^\circ=\Gamma\backslash G/K_0.
		$$
		Here $S$ imposes determinant one on the block diagonal matrices.
		Since $K_0$ is compact, $\operatorname{pr}:\Gamma\backslash G\to Z^\circ$ is proper.
		\begin{step}{Spray leaves expressed as projected $J$-orbits.}
		\end{step}
		Let $\Delta_0\subset\B^n$ be the complex geodesic disk through $o$
		tangent to $\ell_0$, and let $P_\Delta\subset\C^{n+1}$ be its
		defining two-plane, so that
		$$
		\Delta_0=\B^n\cap\PP(P_\Delta),\quad T_o\Delta_0=\ell_0.
		$$
		Let $J\simeq\SU(1,1)$ be a subgroup of $G$ acting on $P_\Delta$ and trivially on
		$P_\Delta^{\perp_h}$. Since $J$ acts transitively on $\Delta_0$,
		its orbit through $(o,\ell_0)$ is precisely the tangent lift
		$$
		J\cdot(o,\ell_0)
		=\{(y,[T_y\Delta_0]):y\in\Delta_0\}
		\subset\PP(T_{\B^n}).
		$$
		For $\Gamma a\in\Gamma\backslash G$, put
		$\zeta_a:=\operatorname{pr}(\Gamma a)\in Z^\circ$.
		The complex geodesic determined by $\zeta_a$ is the image of
		$a\Delta_0$ in $X$, and its tangent lift is
		$\widehat\Delta_{\zeta_a}=\operatorname{pr}(\Gamma aJ)$.
		Since $K_0$ normalizes $J$, this description does not depend on
		the choice of representative of $\zeta_a$. Thus  $\operatorname{pr}(\Gamma aJ)$ is exactly the leaf of the geodesic spray through $\Gamma a K_0$. As $I$ is spray-invariant, $\operatorname{pr}^{-1}(I)$ is closed
		and right $J$-invariant.
		
		\begin{step}{Ratner orbit closures and totally geodesic incidences.}
		\end{step}
		Since $J$ is generated by one-parameter unipotent subgroups,
		for each $\Gamma a\in\operatorname{pr}^{-1}(I)$, Ratner's
		orbit-closure theorem \cite{Ratner} (see also
		\cite[Theorem 1.3]{Shah}) gives a closed subgroup
		$J\subset H\subset G$ with
		$$
		\overline{\Gamma aJ}=\Gamma aH\subset\operatorname{pr}^{-1}(I).
		$$
		The quotient $(a^{-1}\Gamma a\cap H)\backslash H$ has finite Haar volume.
		
		By \cite[Lemma 2.7(2),(4)]{BFMS}, $H\cdot o$ is a totally geodesic $k$-ball, since it contains $J\cdot o=\Delta_0$.
		The group $H$ contains the special unitary group acting on the
		$(k+1)$-plane defining $H\cdot o$ and trivially on its
		orthogonal complement. For $N:=\operatorname{Stab}_G(H\cdot o)$,
		the space $N/H$ is compact and $H$ is transitive on
		$\PP(T_{H\cdot o})$.
		The orbit $\Gamma aN$ is closed \cite[Lemma 8.2(2)]{BFMS}.
		For $\Lambda:=a^{-1}\Gamma a\cap N$, the map
		$\Lambda\backslash N\to\Gamma aN$, $\Lambda n\mapsto\Gamma an$,
		is a homeomorphism \cite[proof of Lemma 2.2]{Shah}.
		With $S:=\operatorname{image}(a(H\cdot o)\to X)$, we obtain
		$$
		\PP(T_{Y_S})\simeq\Lambda\backslash N/(N\cap K_0),
		\quad A_S=\operatorname{pr}(\Gamma aN)
		=\operatorname{pr}(\Gamma aH)\subset I.
		$$
		Quotienting the proper map $\operatorname{pr}|_{\Gamma aN}$
		by the compact group $N\cap K_0$ proves properness of
		the incidence map. Thus $A_S$ is closed in $Z^\circ$.
		Replacing $K_0$ by $K$ proves properness of $\nu_S$;
		\cite[Lemma 8.2(3)]{BFMS} gives finite volume.
		If $k=n$, then $A_S=Z^\circ$, which is excluded.
		
		In summary, the construction starts from the chosen point
		$\zeta_a\in I$ and the spray leaf through it:
		$$
		\zeta_a\in\widehat\Delta_{\zeta_a}
		=\operatorname{pr}(\Gamma aJ)
		\subset A_S=\operatorname{pr}(\overline{\Gamma aJ})\subset I.
		$$
		Moreover, $a(H\cdot o)$ is a totally geodesic
		ball, and $\nu_S:Y_S\to X$ is a proper holomorphic immersion;
		hence $S$ is indeed a totally geodesic submanifold.
		Since $\zeta_a$ was arbitrary, every point of $I$ lies in
		some such $A_S\subset I$.
		
		\begin{step}{Countability and Baire's theorem}
		\end{step}
		The algebra $\operatorname{span}_{\C}\rho_0(\Gamma_{\widetilde S})
		\subset\operatorname{End}(\C^{n+1})$ is spanned by a finite tuple
		of group matrices. Let $P_0$ be the defining plane of $\widetilde S$.
		By Step \ref{step 2222} in the proof of
		Proposition~\ref{prop:wronskian}, $P_0$ is irreducible and
		nonunitary, while $P_0^{\perp_h}$ is unitary.
		Any irreducible submodule on which $h$ is indefinite
		projects nontrivially to $P_0$, hence isomorphically onto it.
		Its projection to $P_0^{\perp_h}$ is zero, since $P_0$ is
		nonunitary. Thus $P_0$ is the unique such submodule, and
		the algebra determines $\widetilde S$. 
		Since $\Gamma$ has countably many finite tuples, there are
		countably many incidences. Hence
		$$
		I=\bigcup_{\alpha\in\mathbb N}A_{S_\alpha},
		\quad A_{S_\alpha}\subset I\text{ closed algebraic}.
		$$
		Baire's theorem gives one member with nonempty relative interior.
		Irreducibility then gives $I=A_{S_\alpha}$.
	\end{proof}

	\begin{remark}
		\label{rem:maximal-invariant-incidences}
		For every $S\in\Tmax(X)$, the tangent incidence $A_S$
		is \emph{maximal} among proper irreducible invariant subvarieties of $Z^\circ$. Indeed, suppose
		$A_S\subset I\subsetneq Z^\circ$, where $I$ is closed, 
		irreducible and invariant.
		By Proposition \ref{prop:invariant-incidence}, we have
		$I=A_T$ for some proper totally geodesic submanifold
		$T\subset X$. Thus the inclusion $A_S\subset I=A_T$ implies $S\subset T$.
		The maximality of $S$ then gives $T=S$, and hence
		$I=A_S$.
	\end{remark}

	\section{Finite descent and degree estimates}
	\label{sec:descent}
	
	In Subsection \ref{subsec:finitedescent}, we start from the
	zero divisor of $p'$ supplied by
	Corollary \ref{cor:initial-divisor} and construct a finite
	descent using tangency to the geodesic spray.
	At each nonterminal step, dimension decreases  without
	losing any invariant subvariety in $Z^\circ$.
	By Proposition \ref{prop:invariant-incidence} and Remark  \ref{rem:maximal-invariant-incidences}, this captures
	the Zariski closure $\overline{A_S}$ of the tangent incidence
	of every maximal totally geodesic submanifold $S$; in particular, the finiteness result  due to \cite{BFMS,BU} (Theorem \ref{thm:known-finiteness}) is recovered. 
	
	In Subsection \ref{subsec:degree-descent}, we bound the total
	$L$-degree of these closures.
	We realize the descent in families with smooth general
	fibers, ensuring that the required components occur with
	positive multiplicity in the special fibers.
	Positive specialization and intersection theory then give
	a degree bound depending only on $n$ and $\deg_L(Z)$.

	\subsection{Finite descent procedure}
	\label{subsec:finitedescent}
	Fix an integer $d$ with $1\le d\le 2n-2$, and let
	\begin{equation}\label{eq:grassmann}
		\mathsf p\colon \mathscr H_d := \operatorname{Gr}_d(T_Z) \longrightarrow Z
	\end{equation}
	be the relative Grassmann bundle of $d$-planes in $T_Z$,
	with universal exact sequence
	\begin{equation}\label{eq:uniquotient}
		0\longrightarrow\mathcal F
		\longrightarrow\mathsf p^*T_Z
		\longrightarrow\mathcal N
		\longrightarrow0.
	\end{equation}
	Here $\mathcal F$ is the \emph{tautological subbundle}
	of rank $d$, with fiber
	$\mathcal F_{(x,F_x)}=F_x\subset T_xZ$, and
	$\mathcal N$ is the \textit{universal quotient bundle}
	of rank $2n-1-d$.
	\begin{definition}\label{def:normalspray}
		Composing the spray $\iota$ defined in Proposition \ref{prop:logarithmic-spray} with the
		universal quotient gives the \emph{normal spray}
		\begin{equation}\label{eq:normal-spray}
			\iota^{\perp}\colon\mathsf p^{*}\OO_Z(-1)
			\xrightarrow{\;\mathsf p^{*}\iota\;}
			\mathsf p^{*}T_{Z}
			\longrightarrow
			\mathcal N.
		\end{equation}
	\end{definition}
	Viewed as a section of $\mathcal N\otimes\mathsf p^*\OO_Z(1)$, the normal spray $\iota^{\perp}$ is regular on $\mathscr H_d$, including over the boundary by Proposition \ref{prop:logarithmic-spray}.
	
	Denote the zero scheme of $\iota^{\perp}$ by $\mathrm V(\iota^{\perp})\subset\mathscr H_d$.
	\footnote{Henceforth, we use $\mathrm V(s)$ similarly for the zero scheme
		of a vector-bundle section $s$.}
	
	\begin{definition}[Gauss map and Nash transform; cf.  for example
		{\cite[pp. 297--298]{Nobile} and 
			\cite[\S2.2]{AluffiGoresky}}]\label{def:gauss-nash}
		Let $W\subset Z$ be an irreducible closed subvariety
		of dimension $d$ with $W\cap Z^\circ\ne\varnothing$.
		\begin{enumerate}
			\item
			When $W$ is smooth, its \emph{Gauss map} is
			$$
			\gamma_W\colon W\longrightarrow\mathscr H_d,\quad
			w\longmapsto(w,T_wW),
			$$
			where $T_wW\subset T_wZ$ denotes the complex
			tangent space of $W$ at $w$.
			
			\item
			In general, $\gamma_W$ is defined on $W_{\reg}$.
			The \emph{Nash transform} of $W$ is
			$$
			\operatorname{Nash}(W)
			:=\overline{\gamma_W(W_{\reg})}^{\,\mathrm{Zar}}
			\subset\mathscr H_d.
			$$
			Then $\operatorname{Nash}(W)$ is irreducible
			of dimension $d$, and $\mathsf p$ restricts to
			a proper birational morphism
			$\operatorname{Nash}(W)\rightarrow W$ that is an
			isomorphism over $W_{\reg}$.
		\end{enumerate}
	\end{definition}
	The restriction of $\iota^\perp$ to $\Nash(W)$ vanishes
	identically if and only if $W\cap Z^\circ$ is invariant.
	The following lemma ensures that invariant subvarieties
	contained in the singular locus of $W$ are also detected.
	
	\begin{lemma}\label{lem:whitney-detection}
		Let $I\subset W\cap Z^\circ$ be an irreducible closed
		subvariety invariant under the spray.
		There exists a dense smooth open subset $I^\circ\subset I$
		such that, for every $x\in I^\circ$ and
		$(x,F_x)\in\Nash(W)$, we have
		$$
		T_xI\subset F_x.
		$$
		Consequently,
		$$
		\mathsf p^{-1}(I^\circ)\cap\Nash(W)
		\subset\mathrm V(\iota^\perp).
		$$
	\end{lemma}
	
	\begin{proof}
		First suppose that $I\not\subset W_{\mathrm{sing}}$.
		Set $I^\circ:=I_{\reg}\cap W_{\reg}$.
		For every $x\in I^\circ$, the unique point of
		$\Nash(W)$ over $x$ is $(x,T_xW)$, and
		$T_xI\subset T_xW$.
		
		Now suppose that $I\subset W_{\mathrm{sing}}$.
		Choose a  Whitney stratification of $W$
		adapted to $I$ \cite{Verdier}.
		Let $U\subset W_{\reg}$ be its dense open stratum,
		and let $I^\circ\subset I_{\reg}$ be its dense
		stratum in $I$.
		Fix $x\in I^\circ$ and $(x,F_x)\in\Nash(W)$.
		Since $\gamma_W(U)$ is dense in $\Nash(W)$, we can
		choose a sequence $w_j\in U$ such that
		$$
		w_j\longrightarrow x,
		\quad
		T_{w_j}W\longrightarrow F_x.
		$$
		Whitney's condition $(a)$ for the pair $(U,I^\circ)$
		\cite[\S19]{Whitney65} gives
		$$
		T_xI=T_xI^\circ\subset F_x.
		$$
		
		In either case, invariance of $I$ gives
		$$
		\iota\bigl(\OO_Z(-1)_x\bigr)\subset T_xI\subset F_x.
		$$
		Thus the image of $\OO_Z(-1)_x$ in $T_xZ/F_x$
		vanishes. By \eqref{eq:normal-spray},
		$\iota^\perp(x,F_x)=0$, proving the assertion.
	\end{proof}

	\begin{proposition}\label{prop:termination}
		Starting with the irreducible components of $\mathrm V(p')$ that meet $Z^\circ$, iterate the following steps.
		\begin{enumerate}
			\item For any current variety $W$, if $W\cap Z^\circ$ is not invariant, replace it by the $\mathsf p(B)$, where $B$ runs through all irreducible components of $\mathrm V(\iota^{\perp})\cap\Nash(W)$, discarding those with $\mathsf p(B)\subset Z\setminus Z^\circ$. We also discard zero-dimensional varieties: by Proposition
			\ref{prop:logarithmic-spray}, none produces a totally geodesic subvariety.
			\item Stop when $W\cap Z^\circ$ is invariant
			for every retained variety $W$.
		\end{enumerate}
		The process terminates after finitely many steps and captures all maximal proper irreducible invariant subvarieties of $Z^\circ$. In particular, there are only finitely many such maximal subvarieties.
	\end{proposition}
	
	\begin{proof}
		If $W\cap Z^\circ$ is noninvariant, then $\iota^{\perp}|_{\Nash(W)}\not\equiv0$.
		For every component $B$ of $\mathrm  V(\iota^{\perp})\cap\Nash(W)$, we have
		\[
		\dim\mathsf p(B)\leq\dim B<\dim W.
		\]
		Now let $I\subset W\cap Z^\circ$ be an irreducible
		invariant subvariety.
		Lemma \ref{lem:whitney-detection} gives
		$$
		I^\circ\subset\bigcup_B\mathsf p(B).
		$$
		Since $\mathsf p$ is proper, the images $\mathsf p(B)$ are
		closed, so the same inclusion holds for $I$.
		The union is finite, so irreducibility gives $I\subset\mathsf p(B)$ for
		at least one component. Thus a target is never lost, including when it
		is contained in $W_{\mathrm{sing}}$.
		
		The dimension strictly decreases at each nonterminal step,
		so the process terminates after finitely many steps.
		By Corollary~\ref{cor:initial-divisor} and
		Proposition~\ref{prop:invariant-incidence}, each maximal $I$
		lies in $\mathrm V(p')$ and is therefore contained in a
		terminal component $W$.
		Maximality gives $W\cap Z^\circ=I$, hence $W=\overline{I}$.
	\end{proof}
	
	\begin{remark}[Finiteness]
		\label{rem:finiteness-totally-geodesic}
		For every $S\in\Tmax(X)$,
		Remark \ref{rem:maximal-invariant-incidences} and
		Proposition \ref{prop:termination} give a terminal
		variety $W\subset Z$ satisfying
		$$
		W\cap Z^\circ=A_S,
		\quad
		W=\overline{A_S}.
		$$
		Since there are only finitely many terminal varieties
		and $\pi(A_S)=S$, we obtain
		$$
		\#\Tmax(X)<\infty.
		$$
		By Remark \ref{rem:finitecover}(a), we recover the finiteness Theorem \ref{thm:known-finiteness} of maximal totally geodesic submanifolds.
	\end{remark}

	\subsection{Specialization and total degree bounds}
	\label{subsec:degree-descent}
	The descent in Subsection \ref{subsec:finitedescent}
	proves finiteness but does not control the total
	$L$-degree of the closures of maximal proper
	invariant subvarieties of $Z^\circ$. 
	In this subsection, we will bound this total degree by realizing the required
	loci as components of special fibers with positive
	multiplicity. The smooth general fibers have controlled
	normal bundles, which allow us to estimate their degrees.

	For a closed subscheme $Q$ of a projective variety, put
	$$
	[Q]_d
	:=\sum_{\dim Q_j=d}
	\operatorname{length}(\OO_{Q,\eta_j})[Q_j],
	$$
	where the $Q_j$ are the irreducible components of $Q$
	with generic points $\eta_j$;
	see \cite[\S1.5]{Fulton}.
	If $Q$ is pure-dimensional of dimension $d$, write
	$[Q]:=[Q]_d$.
	We recall specialization of cycles
	\cite[Tag 0H4H]{Stacks}.
	
	\begin{definition}[Specialization]
		\label{def:specialization}
		Let $Y$ be a projective variety and $(T,0)$ a smooth
		pointed integral curve with generic point $\eta$.
		Let $\alpha_s$ be an algebraic family of effective
		$d$-cycles on $Y$ over $T\setminus\{0\}$, with generic
		cycle
		$$
		\alpha_\eta=\sum_j n_j[V_{j,\eta}].
		$$
		For each $j$, let $\mathscr V_j\subset Y\times T$
		be the scheme-theoretic closure of $V_{j,\eta}$.
		The \textit{specialization} of $\alpha_s$ at $0$
		is the effective $d$-cycle on $Y$ defined by
		$$
		\specialize_0(\alpha_s)
		:=\sum_j n_j[(\mathscr V_j)_0]_d,
		$$
		where $(\mathscr V_j)_0$ denotes the scheme-theoretic
		fiber over $0$.
	\end{definition}
	
	We now introduce families whose smooth fibers have
	controlled normal bundles and whose special fibers
	contain the loci arising in the descent.
	\begin{definition}\label{def:witness}
		A \emph{deformation witness} of dimension $d$ comprises the following data.
		\begin{enumerate}
			\item\label{item:family}
			A proper flat one-parameter family $\mathscr W\subset Z\times T$
			over a smooth pointed algebraic curve $(T,0)$.
			\item\label{item:marked}
			Finitely many \emph{marked} reduced components of the central
			$d$-cycle $[\mathscr W_{0}]_{d}$, each occurring with positive
			integral multiplicity.
		\end{enumerate}
		These data are required to satisfy the following conditions.
		\begin{enumerate}
			\item[(i)]\label{item:smooth}
			The general fiber $W_{t}$ is smooth and embedded in $Z$; it may
			be disconnected.
			\item[(ii)]\label{item:normal}
			For each connected component $V$ of $W_t$,
			the normal bundle $N_{V/Z}$ admits a filtration
			whose graded pieces are
			$$
			L^{\ell_1}|_V,\ldots,L^{\ell_c}|_V,
			\quad c=2n-1-d,\quad \ell_i>0,
			$$
			with the same exponents $\ell_i$ for every component.
			\item[(iii)]\label{item:numerical}
			$[W_{t}]\numeq
			\bigl(\prod_{i=1}^{c}\ell_{i}\bigr)\,c_{1}(L)^{c}\cap[Z]$.
			\footnote{Here $\numeq$ denotes numerical equivalence of cycles
				on $Z$.}  Consequently,
			\[
			\deg_L W_t
			=\left(\prod_{i=1}^{c}\ell_i\right)\deg_L(Z).
			\]
			
		\end{enumerate}
		We call $\ell=(\ell_{1},\ldots,\ell_{c})$ the \emph{normal list} of the
		witness.
	\end{definition}
	
	The initial witness used in the proof is given by the following example.

	\begin{example}[The initial witness]\label{eg:initial}
		Let $p'$ be the section fixed after
		Corollary \ref{cor:initial-divisor}, and choose a general section
		$p'_1\in H^0(Z,L^{\cn{5}})$ independent of $p'$.
		On a Zariski open neighborhood $T\subset\C$ of $0$, consider the pencil
		\[
		\mathscr W^{\mathrm{in}}=\mathrm V(p'+t p'_1)\subset Z\times T,
		\quad W_t^{\mathrm{in}}=\mathrm V(p'+t p'_1).
		\]
		Each fiber is an effective Cartier divisor; hence this relative
		Cartier divisor is flat over $T$. Since $L$ is very ample by
		Lemma \ref{lem:uniform-generation}, Bertini gives a smooth general
		fiber of dimension $d=2n-2$.  After shrinking $T$ around $0$, we may assume that
		$W_t^{\mathrm{in}}$ is smooth for every
		$t\in T\setminus\{0\}$.  For such a fiber,
		\[
		N_{W_t^{\mathrm{in}}/Z}\simeq L^{\cn{5}}|_{W_t^{\mathrm{in}}},
		\quad [W_t^{\mathrm{in}}]\numeq\cn{5}c_1(L)\cap[Z].
		\]
		The special fiber is $W_0^{\mathrm{in}}=\mathrm V(p')$, and its fundamental
		cycle has the form
		\[
		[W_0^{\mathrm{in}}]_d=\sum_\alpha a_\alpha[W_\alpha],
		\quad a_\alpha\in\mathbb Z_{>0}.
		\]
		This gives a deformation witness with normal list
		$(\cn{5})$, whose marked components are the $W_\alpha$ meeting $Z^\circ$. Its total
		degree is
		\[
		\deg_L[W_0^{\mathrm{in}}]_d
		=\sum_\alpha a_\alpha\deg_L(W_\alpha)
		=\cn{5}\deg_L(Z).
		\]
	\end{example}
	The following lemma supplies the positivity needed to
	retain the required components under specialization.
	\begin{lemma}[Positive specialization]\label{lem:positive-specialization}
		\textup{(i)} Let $W$ be a marked component of a
		$d$-dimensional deformation witness $\mathscr W$, and let $\mathscr G$ be the
		closure of the relative Gauss graphs of its smooth fibers in
		$\mathscr H_d\times T$. Then there is $a\in\mathbb Z_{>0}$ such that
		\[
		[\mathscr G_0]_d\geq a[\Nash(W)].
		\]
		Here $\alpha\geq\beta$ means that
		$\alpha-\beta$ is effective.
		
		\smallskip
		\noindent\textup{(ii)}
		Let $\alpha_t$ be a degeneration of effective
		$d$-cycles on a projective variety $Y$ over a smooth
		pointed curve $(T,0)$, and set
		$\alpha_0:=\specialize_0(\alpha_t)$.
		Let $E_1,\ldots,E_h$ be basepoint-free ample line
		bundles, with nonzero sections
		$f_i\in H^0(Y,E_i)$.
		Fix finitely many distinct $(d-h)$-dimensional
		irreducible components $B_j$ of intersections
		$$
		V\cap\bigcap_{i=1}^h\mathrm V(f_i),
		\quad
		V\text{ an irreducible component of }\alpha_0.
		$$
		For general $u_i\in H^0(Y,E_i)$, there are a smooth
		pointed curve $(C,0)$ and a morphism
		$$
		(t,\varepsilon):(C,0)\longrightarrow
		(T\times\C,(0,0)),
		$$
		with $t(s)\ne0$ and $\varepsilon(s)\ne0$ for general
		$s\ne0$, such that the successive intersections
		$$
		\zeta_s
		=\alpha_{t(s)}\cdot
		\mathrm V(f_1+\varepsilon(s)u_1)\cdots
		\mathrm V(f_h+\varepsilon(s)u_h)
		$$
		are proper and satisfy
		\begin{equation}\label{eq:effective-specialization}
			\specialize_0(\zeta_s)\geq\sum_j[B_j].
		\end{equation}
	\end{lemma}
	
	\begin{proof}
		\textup{(i)} There is a dense smooth open subset
		$W^o\subset W$, disjoint from the other central components,
		such that every limiting tangent $d$-plane of the nearby
		smooth fibers over $x\in W^o$ contains $T_xW$
		\cite[\S5, Corollary 1]{Hironaka}.
		Since both have dimension $d$, they coincide.
		Flatness and properness ensure the existence of such
		limits over $W^o$. Thus $\Nash(W)$ is a $d$-dimensional
		component of $\mathscr G_0$.
		As the scheme-theoretic closure of the generic Gauss graph,
		$\mathscr G$ is flat over the smooth curve $T$; hence
		$[\mathscr G_0]_d=\specialize_0([\mathscr G_t])$
		\cite[Tag 0H4J]{Stacks}. Positivity of
		specialization \cite[\S10.1, Example 10.1.8]{Fulton}
		therefore gives its positive integral multiplicity.

		\textup{(ii)} For general $u_i\in H^0(Y,E_i)$,
		the intersections
		$$
		\gamma_\varepsilon
		:=\alpha_0\cdot
		\mathrm V(f_1+\varepsilon u_1)\cdots
		\mathrm V(f_h+\varepsilon u_h)
		$$
		are proper for general $\varepsilon\ne0$.
		At the generic point of each $B_j$, the $h$
		central equations cut the chosen component $V$ in the expected
		codimension. Positivity of local dynamic intersection multiplicities
		\cite[Example 7.1.10(a), Proposition 11.3 \& Example 11.3.1]{Fulton},
		applied with the same general $u_i$ to all $B_j$, gives
		$$\beta_0:=\specialize_{\varepsilon=0}(\gamma_\varepsilon)
		\geq\sum_j[B_j].
		$$
		
		Fix a projective embedding of $Y$. On a dense open
		$U\subset T\times\C$ with $t\ne0$ and $\varepsilon\ne0$,
		the proper intersections of $\alpha_t$ with
		the divisors $\mathrm V(f_i+\varepsilon u_i)$
		define a morphism to a projective Chow variety
		parametrizing effective $(d-h)$-cycles of their common degree.
		Let $\mathscr R$ be the closure of its graph over $T\times\C$.
		Compatibility with $t$-specialization gives
		$(0,\varepsilon,\gamma_\varepsilon)\in\mathscr R$ for general $\varepsilon\ne0$,
		hence $(0,0,\beta_0)\in\mathscr R$.
		Normalize an integral curve in $\mathscr R$ through this point
		whose general point lies over $U$ \cite[Chapter I, \S3]{Kollar};
		choosing a point above $(0,0,\beta_0)$ gives $(C,0)$ and
		$\specialize_0(\zeta_s)=\beta_0$.
		Here $U$ may be replaced by any prescribed dense open subset.
		This proves \eqref{eq:effective-specialization}.
	\end{proof}

	We now use Lemma \ref{lem:positive-specialization}
	to bound the total degree of the closures of maximal
	invariant subvarieties captured by the descent in
	Proposition \ref{prop:termination}.

	\begin{proposition}[Total degree bound]\label{prop:degree-bound}
		Any finite set $\mathscr I$ of closures of maximal
		proper invariant subvarieties of $Z^\circ$ contained
		in $\mathrm V(p')$ satisfies
		\begin{equation}\label{eq:total-degree}
			\sum_{I\in\mathscr I}\degree{I}
			\leq\cn{7}\deg_L(Z).
		\end{equation}
		One may take
		\[
		\cn{6}=(2n)^{2n-1}\bigl(\cn{5}+\cn{3}(2n-1)+1\bigr).
		\]
		\begin{equation}\label{eq:tree-bound}
			\cn{7}=\bigl(1+(2n-1)^2\bigr)^{2n}\cn{6}^{\,2n-1}.
		\end{equation}
	\end{proposition}
	
	\begin{proof}We implement Proposition \ref{prop:termination}
		using deformation witnesses
		(Definition \ref{def:witness}).
		Calling the members of $\mathscr I$ \textit{targets},
		we retain a marked component containing each
		unrecorded target and discard components contained
		in $Z\setminus Z^\circ$.
		Here we remark that discarding a component removes only its marking, while leaving the flat family and its fundamental cycle intact.
		
		We picture the witnesses as a rooted tree:
		\begin{itemize}
			\item
			\textit{Root}: the initial pencil in
			Example \ref{eg:initial}.
			
			\item
			\textit{Children}: witnesses constructed
			from a given witness in one step.
			
			\item
			\textit{Depth}: the maximum number of steps
			along a path from the root, bounded by
			the strict decrease in dimension.
		\end{itemize}
		
		We now complete the proof in the following five steps.
		\setcounter{step}{0}
		\begin{step}{Very ampleness and global generation.}\end{step}\nopagebreak[4]
		Fix a witness of dimension $d$ with normal list
		$\ell=(\ell_1,\ldots,\ell_c)$, where $$c=2n-1-d.$$
		On the Grassmann bundle $\mathsf p:\mathscr H_d\to Z$ of
		\eqref{eq:grassmann}, with universal quotient bundle $\mathcal N$
		of rank $c$ defined in \eqref{eq:uniquotient}, put
		\[
		\mathcal N'=\mathcal N\otimes\mathsf p^*L^{\cn{3}},\quad
		E_d=\det\mathcal N'\otimes\mathsf p^*L^{\cn{3}+1}.
		\]
		By Lemma \ref{lem:uniform-generation}, the bundle
		$T_Z\otimes L^{\cn{3}}$ is globally generated. Its pullback
		surjects onto $\mathcal N'$, so $\mathcal N'$ is globally generated.
		To see that $E_d$ is very ample, let
		$U=H^0(Z,T_Z\otimes L^{\cn{3}})$. The induced quotient
		\[
		U\otimes\OO_{\mathscr H_d}\twoheadrightarrow\mathcal N'
		\]
		defines, by duality, a map
		$\psi:\mathscr H_d\to\Gr_c(U^*)$,
		where the Grassmannian $\Gr_c(U^*)$ parametrizes
		$c$-dimensional subspaces of $U^*$,
		equivalently rank-$c$ quotients of $U$. 
		Its Pl\"ucker line bundle
		$\OO_{\Gr_c(U^*)}(1)$ is very ample
		\cite[Properties 1.29]{Viehweg}, and
		\[
		\psi^*\OO_{\Gr_c(U^*)}(1)=\det\mathcal N'.
		\]
		The map $(\mathsf p,\psi):\mathscr H_d\hookrightarrow
		Z\times\Gr_c(U^*)$ is a closed immersion: over each $x\in Z$
		it identifies $\mathsf p^{-1}(x)$ with the Grassmannian of
		rank-$c$ quotients of $T_xZ\otimes L_x^{\cn{3}}$ inside
		$\Gr_c(U^*)$. As $L^{\cn{3}+1}$ is very ample by
		Lemma \ref{lem:uniform-generation}, the exterior tensor product
		$L^{\cn{3}+1}\boxtimes\OO_{\Gr_c(U^*)}(1)$ is very ample.
		Its pullback by $(\mathsf p,\psi)$ is $E_d$.
		
		The normal spray \eqref{eq:normal-spray} is a section of
		$\mathcal N\otimes\mathsf p^*\OO_Z(1)$. Using
		$\mathcal N^*\otimes\det\mathcal N\simeq
		\bigwedge^{c-1}\mathcal N$, one obtains
		$$
		\begin{aligned}
			\mathcal N^*\otimes\mathsf p^*\OO_Z(-1)\otimes E_d
			&\simeq\bigwedge^{c-1}\mathcal N'\otimes
			\mathsf p^*\bigl(L^{2\cn{3}}\otimes(L\otimes\OO_Z(-1))\bigr).
		\end{aligned}
		$$
		The right-hand side is globally generated: the first factor
		is an exterior power of the globally generated $\mathcal N'$, and
		$L\otimes\OO_Z(-1)=\pi^*\KX^{\cn{2}}$ is globally generated
		by Lemma \ref{lem:uniform-generation}.

		\begin{step}{Proper intersection equations.}\end{step}\nopagebreak[4]
		Choose a finite generating family of global sections
		\begin{equation}\label{global-generated-bndle}
			\{\varphi_\lambda\}\subset
			H^0\!\left(\mathscr H_d,
			\mathcal N^*\otimes\mathsf p^*\OO_Z(-1)\otimes E_d\right)
		\end{equation}
		and contract them with $\iota^{\perp}$:
		\[
		s_\lambda=\langle\varphi_\lambda,\iota^{\perp}\rangle
		\in H^0(\mathscr H_d,E_d).
		\]
		Since the  bundle   in  \eqref{global-generated-bndle}  is globally generated,
		these sections define the same zero scheme as \(\iota^{\perp}\):
		$$
		\mathrm V(\iota^{\perp})
		=\bigcap_\lambda\mathrm V(s_\lambda)
		$$
		scheme-theoretically.
		
		For the subsequent cuts, let \(W\) run through the noninvariant
		marked components, and fix a finite collection of
		\(b\)-dimensional irreducible components
		\[
		B\subset \mathrm V(\iota^{\perp})\cap\Nash(W)
		\]
		which are needed to retain all targets.
		At the generic point of each $B$, the local ring of
		$\Nash(W)$ has dimension $d-b$ and, after trivializing $E_d$,
		the ideal generated by the $s_\lambda$ has radical equal to the maximal ideal.
		Thus, by \cite[Corollary 8.3.9 \& Theorem 8.6.6]{SH}, general linear combinations
		$$
		\sigma_i=\sum_\lambda a_{i\lambda}s_\lambda,
		\quad 1\leq i\leq d-b,
		$$
		generate an ideal with the same radical.
		The finitely many nonempty open conditions admit a common choice
		$a_{i\lambda}\in\C$. Hence each $B$ is a proper intersection component of
		$$
		\Nash(W)\cap \mathrm V(\sigma_1,\ldots,\sigma_{d-b}).
		$$
		
		\begin{step}{Grouping components and constructing a single child witness.}\end{step}\nopagebreak[4]
		Let \(B\subset \mathrm V(\iota^{\perp})\cap\Nash(W)\) be one of the irreducible
		components fixed above, and set \(Y=\mathsf p(B)\) and
		\begin{equation}\label{eq:dimension-pair}
			(b,r_0)=(\dim B,\dim Y).
		\end{equation}
		Since $Y$ contains a positive-dimensional target and
		$\iota^\perp|_{\Nash(W)}\not\equiv0$, one has
		$1\leq r_0\leq b<d$. 
		We group all relevant \(B\) from all noninvariant marked
		components of the fixed witness with the same pair \((b,r_0)\) into a finite
		set \(\mathscr B\) and construct a single child witness for that group.

		Let \(\mathscr G_t\subset\mathscr H_d\) be the Gauss graph of the smooth
		witness fiber \(W_t\). For general
		\(u_i\in H^0(\mathscr H_d,E_d)\), consider
		\[
		\mathscr G_t\cap\bigcap_{i=1}^{d-b}
		\mathrm V(\sigma_i+\varepsilon u_i).
		\]
		
		Since $E_d$ is very ample, successive Bertini intersections
		on $\mathscr G_t\simeq W_t$ are smooth and have pure dimension
		$b$ for general $(t,\varepsilon)$.
		By Lemma \ref{lem:positive-specialization}(i), for each
		marked component $W$, the special Gauss cycle satisfies
		\[
		[\mathscr G_0]_d\geq a_W[\Nash(W)],
		\quad a_W\in\mathbb Z_{>0}.
		\]
		Since each $B\in\mathscr B$ is a proper intersection component of
		$\Nash(W)\cap \mathrm V(\sigma_1,\ldots,\sigma_{d-b})$,
		Lemma \ref{lem:positive-specialization}(ii), applied to
		$\alpha_t=[\mathscr G_t]$ and the divisors $\mathrm V(\sigma_i)$, with the full list $\mathscr B$, gives a smooth pointed
		curve $(T',0)$ as in that lemma.
		Denote these intersections along
		$(t(s),\varepsilon(s))$ by $Q_s$.
		Then there are positive integers
		\(m_B\) such that
		\[
		\specialize_0([Q_s])
		\geq
		\sum_{B\in\mathscr B}m_B[B].
		\]
		We may also require $Q_s$ to be smooth for general $s$.
		Let $\mathscr Q\subset\mathscr H_d\times T'$ be the
		scheme-theoretic closure of these smooth generic fibers.
		As in the proof of Lemma
		\ref{lem:positive-specialization}\textup{(i)}, this closure is flat
		over $\OO_{T',0}$, and specialization gives
		\[
		[\mathscr Q_0]_b=\specialize_0([Q_s]).
		\]
		
		Next choose $b-r_0$ general divisors
		$H_1,\ldots,H_{b-r_0}\in|E_d|$. Successively they avoid the
		associated points of each preceding central fiber and meet the
		smooth general fibers transversely. Their local equations are
		nonzero divisors on the central fibers; the local criterion for
		flatness then preserves flatness after each cut. We obtain
		\[
		\mathscr X:=\mathscr Q\cap\bigcap_{j=1}^{b-r_0}(H_j\times T'),
		\]
		with smooth general fibers $X_s$ of dimension $r_0$ and
		\[
		[\mathscr X_0]_{r_0}
		=[\mathscr Q_0]_b\cdot H_1\cdots H_{b-r_0}
		\geq\sum_{B\in\mathscr B}m_B
		\bigl([B]\cdot H_1\cdots H_{b-r_0}\bigr).
		\]
		Here the intersection cycles are effective, including on the
		special fiber.
		For $Y=\mathsf p(B)$, the generic fiber
		$B_{\eta_Y}$ has dimension $b-r_0$. The restriction of $E_d$ is
		ample, so 
		\[
		e_B:=\int_{B_{\eta_Y}}c_1(E_d)^{b-r_0}>0
		\]
		is the degree of the zero-cycle cut on this generic fiber
		(with $e_B$ the generic degree of
		$\mathsf p|_B:B\to Y$ when $b=r_0$). Therefore
		\begin{equation}\label{eq:positive-projection}
			\mathsf p_*[\mathscr X_0]_{r_0}
			\geq\sum_{Y\in\mathsf p(\mathscr B)}
			\left(\sum_{\substack{B\in\mathscr B\\\mathsf p(B)=Y}}m_Be_B\right)[Y].
		\end{equation}
		
		To obtain a witness in $Z$, use the isomorphism
		$\mathsf p|_{\mathscr G_{t(s)}}:\mathscr G_{t(s)}\simeq W_{t(s)}$.
		It identifies $X_s$ with a smooth embedded complete intersection
		$C_s\subset W_{t(s)}\subset Z$. Let
		$\mathscr C\subset Z\times T'$ be the scheme-theoretic closure
		of the $C_s$ for $s\ne0$. The same argument shows that
		$\mathscr C$ is flat over $\OO_{T',0}$. On the generic fiber,
		$\mathsf p_*[X_s]=[C_s]$; specialization of flat families and its
		compatibility with proper pushforward
		\cite[Tags 0H4J \&  0H4M]{Stacks} yield
		\[
		[\mathscr C_0]_{r_0}
		=\specialize_0([C_s])
		=\specialize_0(\mathsf p_*[X_s])
		=\mathsf p_*[\mathscr X_0]_{r_0}.
		\]
		After shrinking $T'$ around $0$, the families
		$\mathscr Q$, $\mathscr X$, and $\mathscr C$ are flat over $T'$
		and smooth over $T'\setminus\{0\}$.
		By \eqref{eq:positive-projection}, every distinct
		$Y\in\mathsf p(\mathscr B)$ is a marked component of this central
		cycle. Thus one family $\mathscr C$ suffices for the entire group
		$(b,r_0)$. In the next descent step, we form the Gauss graph anew from  the smooth $C_s$.
		
		\begin{step}{Verifying the witness conditions.}\end{step}\nopagebreak[4]
		On the Gauss graph of a smooth $W_t$, the universal quotient
		$\mathcal N$ restricts to $N_{W_t/Z}$. The filtration assumption of
		the witness gives
		$\det N_{W_t/Z}\simeq L^{\sum_i\ell_i}|_{W_t}$.
		Set $\beta:=\sum_i\ell_i+\cn{3}(c+1)+1.$ Then
		\[
		\det\mathcal N'|_{W_t}
		\simeq L^{\sum_i\ell_i+\cn{3}c}|_{W_t},\quad
		E_d|_{W_t}\simeq L^\beta|_{W_t}.
		\]
		All $d-r_0$ divisors used in the previous step therefore restrict to
		$L^\beta$-divisors on $W_t$. Their smooth complete intersection
		$C_s$ has normal bundle
		$N_{C_s/W_t}\simeq(L^\beta|_{C_s})^{\oplus(d-r_0)}$.
		The normal sequence
		\[
		0\longrightarrow(L^\beta|_{C_s})^{\oplus(d-r_0)}
		\longrightarrow N_{C_s/Z}
		\longrightarrow N_{W_t/Z}|_{C_s}\longrightarrow0
		\]
		thus combines with the filtration of $N_{W_t/Z}$ to give the normal list
		\begin{equation}\label{eq:normal-recurrence}
			\ell'=(\underbrace{\beta,\ldots,\beta}_{d-r_0},
			\ell_1,\ldots,\ell_c).
		\end{equation}
		The complete-intersection formula gives
		\[
		[C_s]\numeq\beta^{d-r_0}c_1(L)^{d-r_0}\cap[W_t].
		\]
		Together with flatness and the positively marked components,
		these verify all conditions of Definition \ref{def:witness}.

		\begin{step}{Uniform control.}\end{step}\nopagebreak[4]
		For each $d$, there are at most $(2n-1)^2$ possible pairs
		\eqref{eq:dimension-pair}, regardless of the number of components
		$B$, and every child has strictly smaller dimension. Hence the
		descent tree has depth at most $2n-1$ and at most
		$\bigl(1+(2n-1)^2\bigr)^{2n}$ witnesses.
		Writing $s_\ell=\sum_i\ell_i$ and using $d-r_0\leq2n-1$, \eqref{eq:normal-recurrence} gives
		\[
		s_{\ell'}\leq2n\,s_\ell+
		(2n-1)\bigl(\cn{3}(2n-1)+1\bigr),
		\quad s_{(\cn{5})}=\cn{5}.
		\]
		Hence
		$$
		s_{\ell'}+\cn{3}(2n-1)+1
		\leq2n\bigl(s_\ell+\cn{3}(2n-1)+1\bigr).
		$$
		At depth $q$, this gives
		$$
		s_\ell+\cn{3}(2n-1)+1
		\leq(2n)^q\bigl(\cn{5}+\cn{3}(2n-1)+1\bigr).
		$$
		Since $q\leq2n-1$, every exponent is at most $\cn{6}$;
		each normal list has at most $2n-1$ entries.
		
		By flatness and
		Definition \ref{def:witness}\textup{(iii)}, each witness
		satisfies
		$$
		\deg_L[W_0]_d
		=\deg_L W_t
		=\left(\prod_i\ell_i\right)\deg_L(Z).
		$$ Summing gives
		\[
		\sum_{\text{all witnesses}}\deg_L[W_0]_d
		\leq\bigl(1+(2n-1)^2\bigr)^{2n}
		\cn{6}^{\,2n-1}\deg_L(Z).
		\]
		When a marked component $W$ has invariant interior
		$W\cap Z^\circ$, record it. 
		Maximality identifies that component with any target it contains.
		The target-coverage argument and Proposition \ref{prop:termination}
		ensure that every member of $\mathscr I$ is eventually recorded;
		its positive integral multiplicity ensures that its degree is absorbed into the degree of that witness. Repeated recordings only enlarge the RHS. The displayed estimate proves \eqref{eq:total-degree} with
		\eqref{eq:tree-bound}.
	\end{proof}

	\section{Proof of the main theorem}
	\label{sec:volume}
	\textbf{We return to the original quotient $X$.}
	Let $X'=\Gamma'\backslash\B^n\to X$ be the cover
	chosen in Lemma \ref{lem:deep-cover}, with toroidal
	compactification $(\overline{X'},D')$.
	Use $Z$, $\pi$, and $Z^\circ$ as in \eqref{ZZZ},
	and $L$ as in Lemma \ref{lem:uniform-generation},
	for this cover.
	
	We first compare incidence degrees with volumes on $X'$.
	The total degree bound of
	Proposition \ref{prop:degree-bound} then gives a volume
	estimate on this cover.
	Summing over all lifts cancels the covering degree,
	and the Chern--Gauss--Bonnet formula gives
	Theorem \ref{thm:main}.
	\subsection{Degree bounds the volume}
	
	\begin{lemma}\label{lem:volume-weight}
		Let $S\subset X'$ be a complex totally geodesic submanifold
		of dimension $1\leq k<n$, and let
		$\nu_S:Y_S\to S$ be its normalization. Then
		\begin{equation}\label{eq:degree-volume}
			\degree{\overline{A_S}}\geq\cn{8}\Vol(Y_S),
		\end{equation}
		where
		\[
		\cn{8}:=
		\min_{1\leq k<n}
		\left\{\binom{2k-1}{k-1}
		\left(\frac{\cn{2}}2\right)^k
		k!\left(\frac{n+1}{4\pi}\right)^k\right\}>0.
		\]
	\end{lemma}

	\begin{proof}
		Apply Lemma \ref{lem:toroidal-incidence} to $S\subset X'$,
		and denote the resulting maps by
		$f:\overline{Y_S}\to\overline{X'}$ and $\psi_S$.
		Let
		\[
		\pi_Y:\PP\bigl(T^1_{\overline{Y_S}}(-\log D_{Y_S})\bigr)
		\longrightarrow \overline{Y_S}
		\]
		be the natural projection. Put
		$$
		\xi_Y:=
		c_1\!\left(
		\OO_{\PP(T^1_{\overline{Y_S}}(-\log D_{Y_S}))}(1)
		\right).
		$$
		By Lemma \ref{lem:toroidal-incidence} and
		Proposition \ref{prop:toroidal-properties},
		$\xi_Y=\psi_S^*c_1(\OO_Z(1))$ is nef.
		
		The classes $f^*c_1(K_{\overline{X'}})$ and $f^*c_1\!\left(\OO_{\overline{X'}}(K_{\overline{X'}}+D')\right)$ are nef.
		Moreover, Lemma \ref{lem:deep-cover} gives
		$$
		\begin{aligned}
			&f^*c_1(K_{\overline{X'}})
			-\frac12f^*c_1\!\left(\OO_{\overline{X'}}(K_{\overline{X'}}+D')\right)\\
			&\quad=\frac12f^*c_1\!\left(\OO_{\overline{X'}}(K_{\overline{X'}}-D')\right),
		\end{aligned}
		$$
		which is also nef. Since $\psi_S$ is generically
		of degree one onto $\overline{A_S}$, and
		\[
		\psi_S^*c_1(L)
		=\xi_Y+\cn{2}\pi_Y^*f^*c_1(K_{\overline{X'}}),
		\]
		the nonnegativity of mixed intersections of nef
		classes yields
		\[
		\begin{aligned}
			\degree{\overline{A_S}}
			&=\int_{\PP(T^1_{\overline{Y_S}}(-\log D_{Y_S}))}
			\bigl(\xi_Y+\cn{2}\pi_Y^*f^*c_1(K_{\overline{X'}})\bigr)^{2k-1}\\
			&\geq
			\int_{\PP(T^1_{\overline{Y_S}}(-\log D_{Y_S}))}
			\left(\xi_Y+
			\frac{\cn{2}}2\pi_Y^*f^*c_1\!\left(\OO_{\overline{X'}}(K_{\overline{X'}}+D')\right)\right)^{2k-1}.
		\end{aligned}
		\]
		
		Expanding the last intersection, all terms are
		nonnegative. Keeping the term containing
		$\xi_Y^{k-1}$, we obtain
		\[
		\degree{\overline{A_S}}
		\geq
		\binom{2k-1}{k-1}
		\left(\frac{\cn{2}}2\right)^k\cdot
		\int_{\PP(T^1_{\overline{Y_S}}(-\log D_{Y_S}))}
		\xi_Y^{k-1}\,
		\left(\pi_Y^*f^*c_1\!\left(\OO_{\overline{X'}}(K_{\overline{X'}}+D')\right)\right)^k.
		\]
		The projective-bundle formula
		\cite[\S3.3]{Fulton} gives $\pi_{Y*}(\xi_Y^{k-1})=1.$
		Hence
		\begin{equation}\label{eq:incidence-intersection}
			\degree{\overline{A_S}}
			\geq
			\binom{2k-1}{k-1}
			\left(\frac{\cn{2}}2\right)^k
			\int_{\overline{Y_S}}\left(f^*c_1\!\left(\OO_{\overline{X'}}(K_{\overline{X'}}+D')\right)\right)^k.
		\end{equation}
		
		It remains to identify the last intersection with
		the hyperbolic volume. With our normalization of
		holomorphic sectional curvature $-1$, the first
		Chern form of $K_{X'}$ is
		$(n+1)\omega_{X'}/(4\pi)$.
		Since $f|_{Y_S}=\nu_S$ is a totally geodesic
		isometric immersion, $\nu_S^*\omega_{X'}=\omega_{Y_S}$.
		Furthermore, the boundary compatibility
		\eqref{eq:toroidal-geodesic} allows us to apply
		the logarithmic Chern--Weil formula
		\cite[Proposition 2.2.2]{Cadorel}
		\cite[Theorem 1.4]{Mumford}. Thus
		\begin{equation}\label{eq:restricted-chern-weil}
			\begin{aligned}
				\int_{\overline{Y_S}}\left(f^*c_1\!\left(\OO_{\overline{X'}}(K_{\overline{X'}}+D')\right)\right)^k
				&=\left(\frac{n+1}{4\pi}\right)^k
				\int_{Y_S}\omega_{Y_S}^k\\[0.5em]
				&=k!\left(\frac{n+1}{4\pi}\right)^k
				\Vol(Y_S).
			\end{aligned}
		\end{equation}
		Combining \eqref{eq:incidence-intersection}
		and \eqref{eq:restricted-chern-weil}, and taking
		the minimum over $1\leq k<n$, proves
		\eqref{eq:degree-volume}.
	\end{proof}

	\subsection{Cancellation of the covering degree}
	Put $\mathscr T:=\Tmax(X)$, which is finite by
	Remark \ref{rem:finiteness-totally-geodesic}.
	
	\begin{proposition}\label{prop:cover-volume}
		For every smooth finite-volume non-arithmetic ball quotient
		$X=\Gamma\backslash\B^n$, we have
		\begin{equation}\label{eq:weighted-technical}
			\sum_{S\in\mathscr T}\Vol(Y_S)
			\leq\cn{9}\Vol(X).
		\end{equation}
		Here
		\[
		\cn{9}:=\frac{\cn{7}\cn{4}}{\cn{8}}\,
		n!\left(\frac{n+1}{4\pi}\right)^n.
		\]
	\end{proposition}

	\begin{proof}
		Put $b:=[\Gamma:\Gamma']$.
		Let $\mathscr T'$ be the set of distinct lifts to $X'$
		of the members of $\mathscr T$. Each such lift remains
		maximal, so its tangent incidence is maximal invariant by
		Remark \ref{rem:maximal-invariant-incidences}.
		By Corollary \ref{cor:initial-divisor},
		Proposition \ref{prop:degree-bound} together with
		Lemma \ref{lem:volume-weight}, one has
		\[
		\begin{aligned}
			\cn{8}\sum_{S'\in\mathscr T'}\Vol(Y_{S'})
			&\leq \sum_{S'\in\mathscr T'}\degree{\overline {A_{S'}}}\\[0.5em]
			&\leq \cn{7}\deg_L(Z)\\[0.5em]
			&\leq \cn{7}\cn{4}\int_{\overline{X'}}c_1\!\left(\OO_{\overline{X'}}(K_{\overline{X'}}+D')\right)^n\\[0.5em]
			&=\cn{7}\cn{4}\,n!
			\left(\frac{n+1}{4\pi}\right)^n b\Vol(X).
		\end{aligned}
		\]
		The last inequality is \eqref{eq:ambient-degree},
		and the last equality follows from \eqref{eq:canonical-volume}
		and $\Vol(X')=b\Vol(X)$.
		
		For each $S\in\mathscr T$, fix a lift
		$\widetilde S\subset\B^n$ and put
		$\Gamma_S:=\Gamma_{\widetilde S}$.
		The distinct lifts $S'_\gamma\subset X'$ are indexed by
		$\Gamma'\backslash\Gamma/\Gamma_S$, and the induced maps
		$Y_{S'_\gamma}\to Y_S$ are finite unramified covers of
		degree $[\Gamma_S:\Gamma_S\cap\gamma^{-1}\Gamma'\gamma]$,
		with
		$$
		\sum_{[\gamma]\in\Gamma'\backslash\Gamma/\Gamma_S}
		[\Gamma_S:\Gamma_S\cap\gamma^{-1}\Gamma'\gamma]
		=b.
		$$
		Since hyperbolic volume scales with covering degree,
		it follows that
		\[
		\sum_{S'\mapsto S}\Vol(Y_{S'})
		=b\Vol(Y_S).
		\]
		Summing over $S\in\mathscr T$, we obtain
		\begin{equation}\label{eq:covering-volume}
			\sum_{S'\in\mathscr T'}\Vol(Y_{S'})
			=b\sum_{S\in\mathscr T}\Vol(Y_S).
		\end{equation}
		Substituting \eqref{eq:covering-volume} into the
		preceding degree bound and cancelling $b$ proves
		\eqref{eq:weighted-technical}.
	\end{proof}
	
	\subsection{Proof of Theorem \ref{thm:main}}
	\begin{proof}
		For each $S\in\mathscr T$, put $k:=\dim Y_S$ and let
		$$\chi(Y_S):=\sum_{j=0}^{2k}(-1)^j
		\dim_{\mathbb C}H^j(Y_S,\mathbb C)$$
		denote its topological Euler characteristic.
		The Chern--Gauss--Bonnet formula
		\cite[pp. 719--720]{HP96} gives
		$$
		(-1)^k\chi(Y_S)
		=\frac{(k+1)!}{(4\pi)^k}\Vol(Y_S)>0.
		$$
		Thus $(-1)^k\chi(Y_S)\geq1$, so
		$$
		\Vol(Y_S)\geq\frac{(4\pi)^k}{(k+1)!}\geq\cn{10},
		\quad
		\cn{10}:=\min_{1\leq j<n}\frac{(4\pi)^j}{(j+1)!}>0,
		$$
		since $1\leq k<n$.
		Summing over the finite set $\mathscr T$ and applying
		Proposition \ref{prop:cover-volume}, we obtain
		$$
		\cn{10}\#\mathscr T
		\leq\sum_{S\in\mathscr T}\Vol(Y_S)
		\leq\cn{9}\Vol(X).
		$$
		Setting $C_n:=\cn{9}/\cn{10}$, we conclude that
		$$
		\#\Tmax(X)=\#\mathscr T\leq C_n\Vol(X),
		$$
		as required in \eqref{eq:main-count}.
	\end{proof}

	\appendix
	\renewcommand{\thesection}{\Roman{section}}
	\section{Proofs of Lemmata \ref{lem:deep-cover} and \ref{lem:toroidal-incidence}\label{appendix}}
	\begin{proof}[{Proof of Lemma \ref{lem:deep-cover}}]
		Let $\Lambda\subset\Gamma$ be the subgroup defining the prescribed
		finite cover, and put
		$$
		\Gamma_0:=\bigcap_{\gamma\in\Gamma}\gamma\Lambda\gamma^{-1}.
		$$
		Then $\Gamma_0$ is normal and of finite index in $\Gamma$.
		If $X$ is compact, $K_X$ is ample, and the cover
		$\Gamma_0\backslash\B^n\to X$ satisfies the assertion.
		Assume that $D\ne\varnothing$.
		
		The divisor $M:=\KX+D$ is semiample
		\cite[\S2.2, Main Theorem]{Mok}, hence nef. Adjunction gives
		\begin{equation}\label{eq:toroidal-positivity}
			M|_{D_i}=K_{D_i}\simeq\OO_{D_i}.
		\end{equation}
		Choose a positive rational number $\varepsilon$ such that
		$H:=M-\varepsilon D$ is ample \cite[Theorem A]{BT}.
		We first show that $M-\sum_i\delta_iD_i$ is ample whenever
		$0<\delta_i\leq\varepsilon$ are rational. Put
		$$
		t:=\frac{\min_i\delta_i}{\varepsilon},
		\quad
		R:=(1-t)M-\sum_i(\delta_i-t\varepsilon)D_i.
		$$
		For an irreducible curve $C\not\subset D$, the identity
		$$
		R=(1-t)H+\sum_i(\varepsilon-\delta_i)D_i
		$$
		gives $R\cdot C\geq0$. For $C\subset D_i$, disjointness of
		the boundary components and \eqref{eq:toroidal-positivity} give
		$$
		R\cdot C=-(\delta_i-t\varepsilon)D_i\cdot C\geq0,
		$$
		since $\OO_{D_i}(D_i)$ is anti-ample.
		Thus $R$ is nef, and
		$M-\sum_i\delta_iD_i=tH+R$ is ample.
		
		For each $D_i$, write the central cusp lattice as
		$Z_i=\langle z_i\rangle\simeq\mathbb Z$
		\cite[\S1.1]{Mok}, and choose an integer $N>2/\varepsilon$.
		By residual finiteness \cite[Theorem VII]{Malcev}, we may choose
		a finite-index subgroup $\Gamma'\subset\Gamma_0$
		that is normal in $\Gamma$ and contains none of
		$z_i^r$, $1\leq r\leq N$.
		The induced compactified cover
		$$
		f:(\overline X',D')\longrightarrow(\overline X,D)
		$$
		has ramification indices
		$e_i:=[Z_i:Z_i\cap\Gamma']>N$ over $D_i$.
		The cusp coordinates give
		$$
		f^*M=K_{\overline X'}+D',
		\quad
		f^*D_i=e_i\sum_{D'_j\mapsto D_i}D'_j.
		$$
		Consequently,
		$$
		K_{\overline X'}-D'
		=f^*\left(M-\sum_i\frac{2}{e_i}D_i\right)
		$$
		is ample, since $0<2/e_i<\varepsilon$.
		Moreover, $K_{\overline X'}+D'=f^*M$ is nef, so
		$$
		K_{\overline X'}
		=\frac12\left((K_{\overline X'}-D')
		+(K_{\overline X'}+D')\right)
		$$
		is ample as well. Finally, $\Gamma'\subset\Gamma_0\subset\Lambda$
		gives the required factorization through the prescribed cover.
	\end{proof}
	
	\begin{proof}[{Proof of Lemma \ref{lem:toroidal-incidence}}]
		The cusp groups of $Y_S$ are subgroups of those of $X$,
		hence unipotent by \S\ref{subsubsec:toroidal}.
		The construction recalled there gives a smooth projective
		toroidal compactification $(\overline{Y_S},D_{Y_S})$.
		For $k=\dim Y_S=1$, this is obtained by filling the cusps.
		
		In the Siegel coordinates $(z_1,\ldots,z_n)$ used in
		\eqref{eq:cusp-coordinates}, a Heisenberg translation
		and a unitary change of horizontal coordinates give
		$$
		\widetilde S
		=\B^n\cap\{z_k=\cdots=z_{n-1}=0\}.
		$$
		In these coordinates, $\nu_S$ extends to a holomorphic map
		$$
		f:\overline{Y_S}\longrightarrow\overline X,
		\qquad f^{-1}(X)=Y_S.
		$$
		The central cusp lattices give
		\begin{equation}\label{eq:toroidal-geodesic}
			f^*w_{n,X}=w_{k,Y}^{\,e}u,
			\qquad e\in\mathbb Z_{>0},
		\end{equation}
		where $w_{n,X}$ and $w_{k,Y}$ are local defining functions
		for the boundaries and $u$ is a holomorphic unit.
		The map $f$ is algebraic by \cite[Proposition~15]{GAGA}.
		Its fibers are finite over $X$ because $\nu_S$ is a proper
		immersion, and over $D$ because the horizontal inclusion
		induces immersions of compact abelian varieties.
		Thus the proper morphism $f$ is finite.
		
		By \eqref{eq:toroidal-geodesic},
		$$
		f^*\left(\frac{dw_{n,X}}{w_{n,X}}\right)
		=e\frac{dw_{k,Y}}{w_{k,Y}}+\frac{du}{u}.
		$$
		Using the logarithmic frames \eqref{eq:logvectors},
		along $D_{Y_S}$ we obtain
		$$
		d_{\log}f=
		\begin{pmatrix}
			I_{k-1}&0\\
			0&0\\
			*&e
		\end{pmatrix},
		\qquad
		\rk d_{\log}f=k.
		$$
		Over $Y_S$, the rank is also $k$, since $\nu_S$ is an immersion.
		Hence
		$$
		d_{\log}f:
		T^1_{\overline{Y_S}}(-\log D_{Y_S})
		\hookrightarrow f^*\Tlog
		$$
		is a subbundle injection.
		
		The map $\psi_S$ factors as
		$$
		\PP\bigl(T^1_{\overline{Y_S}}(-\log D_{Y_S})\bigr)
		\hookrightarrow
		\PP(f^*\Tlog)
		\longrightarrow Z.
		$$
		The first arrow is a closed immersion.
		The second is a base change of the finite map $f$.
		Thus $\psi_S$ is finite.
		The factorization also gives the asserted pullback
		identity for $\OO_Z(1)$.
		Since $\nu_S$ has generic degree one and the maps on
		projective tangent fibers are injective,
		$\psi_S$ has generic degree one onto its image.
		
		Moreover,
		$$
		\psi_S^{-1}(Z^\circ)=\PP(T_{Y_S}),
		\qquad
		\psi_S\bigl(\PP(T_{Y_S})\bigr)=A_S.
		$$
		The source is projective and irreducible.
		Its image is therefore $\overline{A_S}$, with
		$$
		\overline{A_S}\cap Z^\circ=A_S,
		\qquad
		\dim\overline{A_S}=2k-1.
		$$
	\end{proof}

	\medskip
	\noindent\textbf{Acknowledgements.}\,\, Both authors would like to thank Professor Kang Zuo for his very helpful comments. They would also like to thank Professor Nicholas Miller for his interest in the first draft. The second author is deeply grateful to  Professors Kwokwai Chan and  Conan Leung, for their constant encouragement and support.

	\medskip
	\noindent\textbf{Declaration on the use of AI.}\,\,
	The first author has been studying the work of
	Baldi--Ullmo \cite{BU} and related papers since 2024,
	while the second author brought the effective algorithmic
	results of Baldi--Urbanik \cite{BUr} into the discussion.
	Inspired by these works and by the finiteness theorem of
	Bader--Fisher--Miller--Stover \cite{BFMS},
	both authors formulated the quantitative finiteness
	problem considered in this paper through their discussions.
	In a series of conversations with us, ChatGPT-6-Astra and
	ChatGPT-6.1-Sol suggested some key ideas for the proof of our main
	theorem, including the introduction of the Wronskian section
	and the use of the geodesic spray.
	These ideas provided the starting point for an
	intersection-theoretic approach to the counting problem.
	The authors then worked together to develop and refine
	the arguments and complete the final proof through
	further discussions with AI.
	AI tools were also used for searching the literature, proving intermediate results
	and polishing the writing.
	All mathematical statements have been independently
	checked by both authors, who take full responsibility
	for the content.
	
	\renewcommand{\bibliofont}{%
		\footnotesize
		\linespread{1.2}\selectfont
	}
	
\end{document}